\documentclass[12pt]{article}
\usepackage{amssymb}
\usepackage{amsmath}
\usepackage{latexsym}
\usepackage{mathabx}
\newcommand{\E}{{\rm I\!E}}

\newcounter{myownsection}
\def\myownsection{\refstepcounter{myownsection} \setcounter{equation}{0}}

\begin{document}
$\;$\\
\begin{center}
{\bf DECOMPOSITIONS OF THE FREE ADDITIVE CONVOLUTION}
\footnote{This work is partially supported by MNiSW research grant P03A 013 30 and by the EU Network QP-Applications,
Contract No. HPRN-CT-2002-00279}\\[50pt]
{\bf Romuald Lenczewski}\\[5pt]
{\it Institute of Mathematics and Computer Science\\
Wroc{\l}aw University of Technology\\
Wybrze\.{z}e Wyspia\'{n}skiego 27\\
50-370 Wroc{\l}aw, Poland\\
e-mail Romuald.Lenczewski@pwr.wroc.pl}\\[30pt]
\end{center}
\begin{abstract}
We introduce and study a new type of convolution of probability measures
called the {\it orthogonal convolution}, which is related to the monotone
convolution. Using this convolution, we derive alternating 
decompositions of the free additive convolution $\mu \boxplus \nu$ of 
compactly supported probability measures in free probabi\-li\-ty.
These decompositions are directly related to alternating decompositions of
the associated subordination functions. In particular, they allow us to compute free additive 
convolutions of compactly supported measures without using free cumulants or R-transforms. 
In simple cases, representations of Cauchy transforms $G_{\mu \,\boxplus\,\nu}(z)$ 
as continued fractions are obtained in a natural way.
Moreover, this approach establishes a clear connection between 
convolutions and products associated with the main notions of independence (free, monotone and
boolean) in noncommutative probability. Finally, our result leads to natural decompositions of 
the free product of rooted graphs.
\\[10pt]
\end{abstract}
\myownsection
\begin{center}
{\sc 1. Introduction}
\end{center}
The addition problem for classically independent random variables leads to the
classical convolution of measures. Namely, if $X_1$ and $X_2$ are independent random variables
with distributions $\mu$ and $\nu$, respectively, then the classical convolution  $\,\mu * \nu\,$
of measures $\mu$ and $\nu$ gives the distribution of $X_1+X_2$.

In free probability there is an analogue of the addition problem which
leads to the free additive convolution of probability measures.
Namely, let $X_1$ and $X_2$ be random variables, i.e. elements of a
noncommutative probability space $({\cal A}, \varphi)$,
where ${\cal A}$ is a unital algebra and $\varphi$ is a linear functional on ${\cal A}$ with
$\varphi(1)=1$, and suppose that $X_1$ and $X_2$ are free with respect to $\varphi$. If $\mu$ and $\nu$
denote the $\varphi$-distributions of $X_1$ and $X_2$, respectively, then the $\varphi$-distribution
of $X_1+X_2$, denoted $\mu \boxplus \nu$, is called the {\it free additive convolution} of $\mu$ and
$\nu$. This convolution was introduced by Voiculescu [24] for compactly supported
probability measures on the real line. In the procedure of computing $\,\mu\boxplus\nu\,$,
a central role is played by the Cauchy transforms of probability measures.
In particular, if $X$ is a bounded self-adjoint operator
on some Hilbert space ${\cal H}$, then
\begin{equation}
G_{\mu}(z)=\int_{-\infty}^{\infty}\frac{\mu(dx)}{z-x}=\sum_{n=0}^{\infty}\mu(X^n)z^{-n-1},
\end{equation}
where $z$ lies in the open upper half plane ${\mathbb C}^{+}$,
is the Cauchy transform of $\mu$ with moments $\mu(X^n)=\varphi(X^{n})$.

Voiculescu introduced the {\it R-transform} of $\mu$ defined by 
$R_{\mu}(z)=G_{\mu}^{-1}(z)-1/z$, where $G_{\mu}^{-1}(z)$ is the right inverse of 
$G_{\mu}(z)$ with respect to the composition of formal power series (with 
coefficients $r_{\mu}(n)$ called the {\it free cumulants} of $\mu$).
On a suitable domain, $R_{\mu}(z)$ becomes a holomorphic function.
If R-transforms are used to compute $\,\mu\boxplus\nu\,$, one has to invert the Cauchy transforms
$G_{\mu}(z)$ and $G_{\nu}(z)$, which gives $R_{\mu}(z)$ and $R_{\nu}(z)$, add these up
to get $R_{\mu\,\boxplus\, \nu}(z)=R_{\mu}(z)+R_{\nu}(z)$ and then invert $G^{-1}_{\mu\,\boxplus \,\nu}(z)$ back to
obtain $G_{\mu\, \boxplus \,\nu}(z)$. Finally, using the Stieltjes inversion formula, one
can compute $\mu \boxplus \nu$ (for details, see [24] and [25]).

The additivity of the R-transform is analogous to the additivity of
the logarithm of the Fourier transform $\,{\cal F}_{\mu}(it)\,$, or of
the associated exponential moment generating function
\begin{equation}
{\cal F}_{\mu}(z)=\sum_{n=0}^{\infty}\frac{\mu(X^n)}{n!}z^{n}
\end{equation}
for the measure $\mu$, and the free cumulants are the analogues of the classical cumulants
which appear in the power series representing $\,{\rm log}\,{\cal F}_{\mu}(it)$.
Nevertheless, in classical probability, one can express the moments of $\mu *\nu$ directly
in terms of the moments of $\mu$ and $\nu$ without using the classical cumulants since
there is a `complete decomposition'
\begin{equation}
{\cal F}_{\mu*\nu}(z)={\cal F}_{\mu}(z)\cdot {\cal F}_{\nu}(z).
\end{equation}
In this paper we find some analogues of the above formula for the convolution $\mu\boxplus \nu$,
which allow us to compute it without using the free cumulants (or the $R$-transform).

For that purpose we shall use the reciprocal Cauchy transforms
of probability measures.
By the {\it reciprocal Cauchy transform} of $\mu$ we understand
\begin{equation}
F_{\mu}(z)=\frac{1}{G_{\mu}(z)}
\end{equation}
and the class of reciprocal Cauchy transforms of Borel probability measures on the real line
${\cal M}$ we denote by ${\cal RC}$.
In fact, they played a central role in the approach 
of Maassen [18] who extended the definition of the additive free
convolution to measures with finite variance. To all measures from class ${\cal M}$, the definition
was later extended by Bercovici and Voiculescu [5]. Another important result in the context
of reciprocal Cauchy transforms is the subordination property, namely that there exist unique functions
$F_{1}, F_{2}\in {\cal RC}$, called {\it subordination functions}, such that
\begin{equation}
F_{\mu \, \boxplus \, \nu}(z)=F_{\mu}(F_{1}(z))=F_{\nu}(F_{2}(z))
\end{equation}
for $z\in{\mathbb C}^{+}$ (proved by Voiculescu [25]
for compactly supported measures and by Biane [6] in the general case).
These functions play a key role in the recent work of 
Belinschi [4], where complex analytic methods are used to give a detailed study of free convolutions.

The above relation was also used by Chistyakov and Goetze [11] who proved that
there exist unique functions $F_{1},F_{2}\in {\cal RC}$ such that
\begin{equation}
F_{\mu}(F_{1}(z))=F_{\nu}(F_{2}(z))\;\;\; {\rm and} \;\;\; F_{\mu}(F_{1}(z))=F_{1}(z)+F_{2}(z)-z,
\end{equation}
for any $\mu, \nu \in {\cal M}$, and therefore one can define $\mu\boxplus \nu$
by any of the equations (1.5). It is worth mentioning here that formulas (1.6) were also 
used by Quenell [21] and Gutkin [12] in their study of free products of graphs and their spectra. 
In their approach, the functions $F_{1}(z)$ and $F_{2}(z)$ correspond to
the (root) spectral distributions of the `branches' of the free product of rooted graphs.

The above formulas are related to the {\it boolean convolution} of measures which appeared in the addition
problem for boolean independent random variables [22] and the {\it monotone convolution}
corresponding to the notion of monotone independence [19,20].
These convolutions can be defined using the reciprocal Cauchy transforms by the equations
\begin{equation}
F_{\mu\vartriangleright \nu}(z)=F_{\mu}(F_{\nu}(z))\;\;\; {\rm and}
\;\;\;
F_{\mu \uplus \nu}(z)=F_{\mu}(z)+F_{\nu}(z)-z.
\end{equation}
where  $\mu \uplus \nu$ and $\mu\vartriangleright \nu$ denote the boolean and monotone
convolutions of $\mu, \nu \in {\cal M}$.
Using these convolutions and (1.5), equations (1.6) correspond to what might be 
viewed as `monotone' and `boolean' decompositions of the free additive convolution. 

Nevertheless, functions $F_{1}(z)$ and $F_{2}(z)$ still depend on both
$\mu$ and $\nu$. In this paper we find and study `complete' decompositions of
the corresponding measures, which express them in terms of a new `basic'
convolution which resembles the monotone convolution.
This new convolution is denoted $\mu \vdash \nu$ and called the 
{\it orthogonal convolution} of $\mu$ and $\nu$. If we continue to use the convenient 
language of transforms, we can define it by its reciprocal Cauchy transform
\begin{equation}
F_{\mu \,\vdash\,\nu}(z)=F_{\mu}(F_{\nu}(z))-F_{\nu}(z)+z.
\end{equation}
which shows how much it differs from the monotone convolution.
Equivalently, the orthogonal convolution of $\mu$ and
$\nu$ can be defined as the unique measure $\,\mu \vdash \nu\,$
determined by the equation
\begin{equation}
K_{\mu\,\vdash\,\nu}(z)=K_{\mu}(F_{\nu}(z))=K_{\mu}(z-K_{\nu}(z)),
\end{equation}
where the $K$-transform of measure $\mu$ is given by $K_{\mu}(z)=z-F_{\mu}(z)$ 
(introduced by Speicher and Woroudi [22]).
The Hilbert-space realization of the orthogonal convolution 
involves projections $P_{1}^{\perp}$ and $P_{2}$ onto 
${\cal H}\ominus {\mathbb C}\xi_{1}$ and ${\mathbb C}\xi_{2}$, respectively, where
$\xi_{1}$ and $\xi_{2}$ are cyclic vectors, which motivates our terminology.
The orthogonal convolution is the main building block of the decompositions of
$\mu \boxplus \nu$ studied in this paper. In particular, we study its
combinatorics and decompose the sum $X_1+X_2$ of free random variables on the
free product of Hilbert spaces as infinite sums of replicas of 
$X_{1}$ and $X_{2}$ which correspond to the `factors' of the `complete' decompositions of $\mu
\boxplus \nu$.

In these decompositions, proven here for compactly supported measures,
the boolean or monotone convolutions deduced from (1.8) are followed by infinite sequences of
orthogonal convolutions of alternating $\mu$ and $\nu$.
Using reciprocal Cauchy transforms and $K$-transforms, the first decomposition
corresponds to the `continued composition' form
\begin{equation}
F_{\mu \,\boxplus \,\nu}(z)=F_{\mu}(z-K_{\nu}(z-K_{\mu}(z-K_{\nu}(z-\ldots ))))
\end{equation}
where the right-hand side is understood as the uniform limit on compact subsets of ${\mathbb C}^{+}$
(a `twin-like' formula is obtained by interchanging $\mu$ and $\nu$).
Viewing the $K$-transform as a slight modification of the reciprocal Cauchy
transform, we can treat this expression as the `monotone-orthogonal decomposition'
of $F_{\mu\, \boxplus \,\nu}(z)$ since it begins with the  composition
of $F$-transforms.
It is not hard to see that this decomposition is closely related to continued fractions.
Another decomposition of $F_{\mu\,\boxplus\,\nu}(z)$ is called the `boolean-orthogonal
decomposition' since it corresponds to the second equation of (1.6). 

Let us remark that the approximants of these decompositions correspond to 
approximations of freeness which were studied from the point of view 
of product states, limit theorems and Gaussian operators. 
Thus, the approximants of the boolean-orthogonal decomposition 
correspond to the {\it hierarchy of freeness} [14,15], whereas
the approximants of the monotone-orthogonal decomposition correspond
to the {\it monotone hierarchy of freeness} [17].
In a different direction goes [16], where a noncommutative extension of the 
Fourier transform was constructed which extends both the Fourier transform and the 
$K$-transform. For other interpolations involving the free
additive convolution, see [7-9]. In turn, the free multiplicative convolution 
and its decompositions will be treated in a separate paper.

The paper is organized as follows. In Section 2 we introduce basic notions.
Some useful combinatorics is developed in Section 3.
In Section 4 we introduce and study the concept of
`orthogonal subalgebras' of a noncommutative probability space and related 
`orthogonal structures' (product of Hilbert spaces, product of $C^{*}$-algebras and 
convolution)
The notion of `orthogonal convolution' is studied in more detail in Section 5 (algebraic properties)
and Section 6 (transforms and analytic properties). Then, in Section 7, we introduce and 
study `subordinate structures' related to subordination functions and 
based on the concept of `s-free subalgebras' of a noncommutative probability space.
We show in Section 8 that they generalize the `branches' of the 
free product of graphs. We also derive an alternating orthogonal decomposition of 
the s-free convolution. This leads to monotone-orthogonal and boolean-orthogonal
decompositions of the free additive convolution given in Section 9.
\\[10pt]
\myownsection
\begin{center}
{\sc 2. Preliminaries}
\end{center}
This section contains preliminaries concerning transforms of probability measures
(more generally, distributions of random variables) and their convolutions
associated with notions of noncommutative independence (free, monotone and boolean).

By a {\it non-commutative probability space} we understand a pair $({\cal A}, \varphi)$,
where ${\cal A}$ is a unital algebra over ${\mathbb C}$ and $\varphi$ is
a linear functional $\varphi: {\cal A}\rightarrow {\mathbb C}$ such that
$\varphi(1)=1$. If ${\cal A}$ is a unital *-algebra and $\varphi$ is positive (called a {\it state}),
then $({\cal A},\varphi)$ is called a *-{\it probability space}. If, in addition,
${\cal A}$ is a $C^{*}$-algebra, then $({\cal A}, \varphi)$ is called a
$C^{*}$-{\it probability space}. By the Gelfand-Naimark-Segal theorem, a $C^{*}$-probability
space can always be realized as a subalgebra of bounded operators on a Hilbert space
${\cal H}$ with a distinguished unit vector $\xi$, for which $\varphi(a)=\langle a\xi, \xi \rangle$
for $a\in {\cal A}$. 

By a {\it random variable} we will understand any element $a$ of the considered algebra ${\cal A}$.
If ${\cal A}$ is equipped with an involution, then a random variable $a$ will be called
{\it self-adjoint} if $a=a^{*}$. The $\varphi$-{\it distribution} of a random variable $a$ is
the functional $\mu_{a}:{\mathbb C}[X]\rightarrow {\mathbb C}$
given by $\mu_{a}(1)=1$, $\mu_{a}(X^{n})=\varphi(a^{n})$. In particular,
if $({\cal A}, \varphi)$ is a $C^{*}$-probability space, then
the distribution $\mu_{a}$ of a self-adjoint random variable $a\in {\cal A}$ extends to
a compactly supported probability measure $\mu$ on the real line.
In that case we will often use the same notation $\mu$ for both the
distribution of $a$ and the associated compactly supported probability measure.

The additive free convolution of distributions (measures)
$\mu \boxplus \nu$ is related to the notions of freeness and free product
of $C^{*}$-algebras  [3,23].
Let $({\cal A}, \varphi)$ be a non-commutative probability space and let
${\cal A}_{i}$, $i\in I$, be unital subalgebras of ${\cal A}$.
The family $({\cal A}_{i})_{i\ \in I}$ is called {\it free} with respect to $\varphi$ if
\begin{equation}
\varphi(a_{1}a_{2}\ldots a_{n})=0
\end{equation}
whenever $a_{j}\in {\cal A}_{i_{j}}\cap {\rm Ker}\varphi$
with $i_{1}\neq i_{2}\neq \ldots \neq i_{n}$. A family of elements $(a_{i})_{i\in I}$
of ${\cal A}$ is called {\it free} if the family of unital subalgebras $({\cal A}_{i})_{i\in I}$
of ${\cal A}$, each generated by $a_{i}$, is free.

In turn, the monotone convolution is related to {\it monotone independence} [19]
which can be defined if the set $I$ is totally ordered. Thus, random variables
$(a_{i})_{i\in I}$ are {\it monotone independent} w.r.t. $\varphi$ if
\begin{equation}
\varphi(a_{i_1}\ldots a_{i_k}\ldots a_{i_n})=
\varphi(a_{i_k})\varphi(a_{i_1}\ldots a_{i_{k-1}}a_{i_{k+1}}\ldots a_{i_n})
\end{equation}
whenever $i_{k-1}<i_{k}$ and $i_{k}>i_{k+1}$, with the understanding that only one of these inequalities
holds if $k\in \{1,n\}$.

In particular, we will say that the pair $(a,b)$ of elements of
${\cal A}$ is {\it monotone independent} w.r.t. $\varphi$ if $a=a_{i}$ and $b=a_{j}$ with $i<j$ and $a_{i}, a_{j}$
are monotone independent w.r.t. $\varphi$. In that case, if the
$\varphi$-distributions of $a$ and $b$ are $\mu$ and $\nu$, respectively, then
the $\varphi$-distribution of $a+b$, denoted $\mu\vartriangleright\nu$, is
called the {\it monotone convolution} of $\mu$ and $\nu$. If $\mu$ and
$\nu$ are probability measures on ${\mathbb R}$, then
$\mu\vartriangleright\nu$ is the unique probability measure on ${\mathbb R}$
which satisfies the equation
\begin{equation}
F_{\mu \vartriangleright
\nu}(z)=F_{\mu}(F_{\nu}(z))
\end{equation}
where $F_{\mu}(z)$ is the reciprocal Cauchy transform of $\mu$.
For details, see [20].

The third convolution, which plays an important role in our approach, is associated
with the so-called {\it boolean independence} [22]. Namely, random variables $(a_{i})_{i\in I}$
are called {\it boolean independent} w.r.t. $\varphi$ if
\begin{equation}
\varphi(a_{i_1}a_{i_2}\ldots a_{i_n})=
\varphi(a_{i_1})\varphi(a_{i_2})\ldots \varphi(a_{i_{n}})
\end{equation}
whenever $i_{1}\neq i_{2}\neq \ldots \neq i_{n}$.

In particular, if two random variables, $a_{1}=a$ and $a_{2}=b$, have $\varphi$-distributions $\mu$ and $\nu$,
respectively, and are boolean independent w.r.t. $\varphi$, then the $\varphi$-distribution of the sum $a+b$
is denoted $\mu\uplus\nu$ and is called the {\it boolean convolution} of $\mu$ and $\nu$.
If $\mu$ and $\nu$ are probability measures on ${\mathbb R}$, then $\mu\uplus\nu$
is the unique probability measure on ${\mathbb R}$ which satisfies the equation
\begin{equation}
K_{\mu \uplus \nu}(z)=K_{\mu}(z)+K_{\nu}(z)
\end{equation}
where $K_{\mu}(z)=z-F_{\mu}(z)$ is the so-called $K$-{\it transform} of $\mu$.
If $\mu$ is a probability measure, then $K_{\mu}:\; {\mathbb C}^{+}\rightarrow {\mathbb C}^{-}\cup {\mathbb R}$
is a holomorphic function, where ${\mathbb C}^{+}$ and ${\mathbb C}^{-}$ denote the open
upper and lower complex half-planes, respectively. For details, see [22].\\[10pt]

\myownsection
\begin{center}
{\sc 3. Combinatorics}
\end{center}
In this section we describe the combinatorics which appears in a natural way in the context
of the orthogonal convolution.

We adopt the following notations.
By ${\cal I}(n)$ we denote the lattice of {\it interval partitions}
of the set $\{1,2, \ldots, n\}$. Thus, any $\pi\in {\cal I}(n)$ is of the form
$\pi=\{\pi_{1}, \pi_{2}, \ldots , \pi_{r}\}$, where
$\pi_{1}\cup\pi_{2}\cup \ldots \cup \pi_{r}=\{1,2, \ldots , n\}$ and
$\pi_{1}<\pi_{2}<\ldots < \pi_{r}$, where blocks are ordered in the natural way, i.e.
$\pi_{i}<\pi_{j}$ means that $k<l$ for all $k\in \pi_{i}$ and $l\in \pi_{j}$.
Note that there is a bijection between ${\cal I}(n)$ and the set of {\it ordered
partitions} of the number $n\in {\mathbb N}$. Namely, the bijection is given by
$\pi\rightarrow (j_{1}, j_{2}, \ldots , j_{r})$, where $j_{i}=|\pi_{i}|$ for $i=1, \ldots, n$.
Clearly, $j_{1}+j_{2}+\ldots +j_{r}=n$.
This bijection will be used in the sequel and both objects, the interval partition $\pi$ and
the corresponding tuple $(j_{1},j_{2}, \ldots , j_{r})$ will be denoted by $\pi$.

For $r$ odd and $\pi\in {\cal I}(n)$ with blocks ordered in the natural way,
we shall use the {\it alternating} decomposition of
$\pi\in {\cal I}(n)$ of the form $\pi=\pi'\cup \pi''$, where
\begin{equation}
\pi'=\{\pi_{1}, \pi_{3}, \ldots , \pi_{r}\},\;\;
{\rm  and}\;\; \pi''=\{\pi_{2},\pi_{4}, \ldots , \pi_{r-1}\},
\end{equation}
with the associated tuples $(j_{1},j_{3}, \ldots , j_{r})$ and $(j_{2},j_{4}, \ldots ,
j_{r-1})$, respectively (a similar definition can be given for $r$ even, but we will not
need it). Further, we write $\pi \preceq \sigma$ for $\pi, \sigma \in {\cal I}(n)$
if $\pi $ is a (not necessarily proper) refinement
of $\pi$. Finally, if $\pi\in {\cal I}(n)$, then by ${\cal I}(\pi)$ we
denote the family of all (interval) subpartitions $\sigma\in {\cal I}(n)$ of the
partition $\pi$ and by ${\cal I}_{{\rm odd}}(\pi)$ - its subset consisting
of such (interval) subpartitions of $\pi$ which are obtained from $\pi$ by
decomposing its every block $\pi_{k}$ into an odd number of subblocks -
these subpartitions will be called {\it odd}.

We will use multiplicative functions on partially ordered sets.
The partially ordered set used here will be the union of lattices $\bigcup_{n\geq 1}{\cal I}(n)$
with the natural partial order, again denoted $\pi\preceq \sigma$, iff
there exists $n\in {\mathbb N}$ such that $\pi, \sigma \in {\cal I}(n)$ and
it holds that $\pi \preceq \sigma$ for $\pi, \sigma$ treated as elements of
${\cal I}(n)$. In particular, for any distribution $\mu$ we define the {\it
moment function}
\begin{equation}
m_{\mu}:\; \bigcup _{n\geq 1} {\cal I}(n)\rightarrow {\mathbb R}
\end{equation}
\begin{equation}
m_{\mu}(\pi)=\mu(j_1)\mu(j_2)\ldots \mu(j_r)
\end{equation} where
$j_{i}=|\pi_{i}|$, $1\leq i\leq r$ and $(\mu(n))_{n\in {\mathbb N}}$ is the
collection of moments of $\mu$. Related to the moment function is the
inverse boolean cumulant function defined below (it differs from the usual
boolean cumulant function $k_{\mu}$ [22] with summation extending over $\pi\preceq \sigma$).\\
\indent{\par}
{\sc Definition 3.1.} Let
the moment function $m_{\mu}$ be related to the multiplicative functions
$k_{\mu}^{*}$ on the lattice $\bigcup_{n\geq 1}{\cal I}(n)$
by the formula
\begin{equation}
m_{\mu}(\sigma)= \sum_{\pi \succeq
\sigma}k_{\mu}^{*}(\pi).
\end{equation}
Then, functions
$k_{\mu}^{*}$, given by the M\"obius inversion formula
\begin{equation}
k_{\mu}^{*}(\pi)=\sum_{\sigma
\succeq \pi}(-1)^{|\sigma|-|\pi|}m_{\mu}(\sigma)
\end{equation}
will be called the {\it inverse boolean
cumulant function}. Using the representation of $\pi$ as
$(j_{1},j_{2}, \ldots , j_{n})$, we can write
$k_{\mu}^{*}(\pi)=k_{\mu}^{*}(j_{1},j_{2},\ldots , j_{n})$.\\
\indent{\par}
{\sc Example 3.1.}
Let us evaluate $k_{\mu}^{*}(\pi)$ for the simplest
partitions $\pi$. We get
\begin{eqnarray*} k_{\mu}^{*}(n)&=&\mu(n)\\
k_{\mu}^{*}(n,m)&=& \mu(n)\mu(m)-\mu(n+m)\\ k_{\mu}^{*}(n,m,k)&=&
\mu(n)\mu(m)\mu(k)-\mu(n+m)\mu(k)- \mu(n)\mu(m+k)+\mu(n+m+k)\\
k_{\mu}^{*}(n,m,k,l)&=&
\mu(n)\mu(m)\mu(k)\mu(l)-\mu(n+m)\mu(k)\mu(l)-\mu(n)\mu(m+k)\mu(l)\\
&&-\mu(n)\mu(m)\mu(k+l) + \mu(n+m)\mu(k+l) +\mu(n+m+k)\mu(l)\\
&&+\mu(n)\mu(m+k+l)-\mu(n+m+k+l).
\end{eqnarray*} where $n,m,k,l\in
{\mathbb N}$. \\
\indent{\par}
A partition $\pi=\{\pi_{1},\pi_{2}, \ldots ,
\pi_{r}\}$ of the set $\{1,2, \ldots , n\}$ is called {\it non-crossing} if
there do not exist numbers $i<k<j<l$ such that $i,j\in \pi_{p}$, $k,l\in
\pi_{q}$ and $p\neq q$. By ${\cal NC}(n)$ we denote the family of
non-crossing partitions of the set $\{1,2, \ldots , n\}$. If $\pi\in {\cal
NC}(n)$, then its block $\pi_{p}$ is {\it inner} with respect to block
$\pi_{q}$ if $i<k<j$ for every $k\in \pi_{p}$ and $i,j\in \pi_{q}$ (then
$\pi_{q}$ is called {\it outer} w.r.t. $\pi_{p}$). Let $o(\pi_{p})$ be the
number of blocks of $\pi$ which are outer w.r.t. $\pi_{p}$. Then the {\it
depth} of $\pi_{p}$ is defined as $d(\pi_{p})=o(\pi_{p})+1$ and
$d(\pi)={\rm max}_{1\leq j \leq r}d(\pi_{j})$ is called the {\it depth} of
$\pi$. By ${\cal NC}_{d}(n)$ we shall denote the family of non-crossing
partitions of depth smaller or equal to $d$. In particular, ${\cal
NC}_{1}(n)={\cal I}(n)$ for every $n\in {\mathbb N}$.

Let us introduce a suitable subfamily of ${\cal NC}_{2}(n)$.\\
\indent{\par}
{\sc Definition 3.2.}
A partition $\pi \in {\cal NC}_{2}(n)$ is called {\it decomposable}
if it can be decomposed as  $\pi=\pi'\cup \pi''$, where
\begin{equation}
\pi'=\{\pi_{1}',\pi_{2}', \ldots , \pi_{p}'\},\;\;\;
\pi''=\{\pi_{1}'', \pi_{2}'', \ldots , \pi_{q}''\}
\end{equation}
with $p\geq 1$ and $q\geq 0$, consisting of $p$ blocks of depth 1,
$\pi_{1}'< \pi_{2}'< \ldots < \pi_{p}'$, and $q$ blocks of depth 2,
$\pi_{1}''<\pi_{2}''< \ldots <\pi_{q}''$, such that 
blocks $\pi_{i}''$ and $\pi_{i+1}''$ are not neighbors.
Denote by ${\cal D}_{2}(n)$ the family of all decomposable partitions 
of the set $\{1,2, \ldots , n\}$.\\
\indent{\par}
{\sc Definition 3.3.}
Let $\pi \in {\cal D}_{2}(n)$ be given with decomposition $\pi=\pi'\cup \pi''$
of the form (3.6). Denote by ${\cal P}(\pi)$ the set of refinements $\eta$ of $\pi$
of the form $\eta=\eta' \cup \eta''$, where $\eta''=\pi''$ and
$$
\eta'=\{\eta_{1}',\eta_{2}', \ldots , \eta_{r}'\}, \;\;\;
$$
is a refinement of $\pi'$ which satisfies the conditions:
(i) $\eta_{1}'<\eta_{2}'<\ldots <\eta_{r}'$, (ii) if two consecutive numbers $i,i+1$
belong to the same block of $\pi'$, they must belong to the same block of $\eta'$.
By a {\it decomposition pair} we understand any pair $(\pi,\eta)$, where
$\pi\in {\cal D}_{2}(n)$ and $\eta\in {\cal P}(\pi)$. Denote by
${\cal DP}_{2}(n)$ the family of decomposition pairs $(\pi, \eta)$, where
$\pi\in {\cal D}_{2}(n)$ and $\eta\in {\cal P}(\pi)$.  \\
\unitlength=1mm
\special{em:linewidth 0.4pt}
\linethickness{0.4pt}
\begin{picture}(90.00,80.00)(-40.00,-25.00)

\put(0.00,40.00){\line(1,0){40.00}}
\put(45.00,40.00){\line(1,0){35.00}}
\put(0.00,30.00){\line(0,1){10.00}}
\put(5.00,30.00){\line(0,1){10.00}}
\put(20.00,30.00){\line(0,1){10.00}}
\put(25.00,30.00){\line(0,1){10.00}}
\put(40.00,30.00){\line(0,1){10.00}}
\put(45.00,30.00){\line(0,1){10.00}}
\put(60.00,30.00){\line(0,1){10.00}}
\put(80.00,30.00){\line(0,1){10.00}}

\put(10.00,30.00){\line(0,1){5.00}}
\put(15.00,30.00){\line(0,1){5.00}}
\put(30.00,30.00){\line(0,1){5.00}}
\put(35.00,30.00){\line(0,1){5.00}}

\put(50.00,30.00){\line(0,1){5.00}}
\put(55.00,30.00){\line(0,1){5.00}}
\put(65.00,30.00){\line(0,1){5.00}}
\put(70.00,30.00){\line(0,1){5.00}}
\put(75.00,30.00){\line(0,1){5.00}}

\put(10.00,35.00){\line(1,0){5.00}}
\put(30.00,35.00){\line(1,0){5.00}}
\put(50.00,35.00){\line(1,0){5.00}}
\put(65.00,35.00){\line(1,0){10.00}}

\put(11.00,27.00){\footnotesize $\pi_{1}''$}
\put(31.00,27.00){\footnotesize $\pi_{2}''$}
\put(51.00,27.00){\footnotesize $\pi_{3}''$}
\put(69.00,27.00){\footnotesize $\pi_{4}''$}
\put(20.00,42.00){\footnotesize $\pi_{1}'$}
\put(62.00,42.00){\footnotesize $\pi_{2}'$}

\put(0.00,15.00){\line(1,0){25.00}}
\put(45.00,15.00){\line(1,0){15.00}}
\put(0.00,5.00){\line(0,1){10.00}}
\put(5.00,5.00){\line(0,1){10.00}}
\put(20.00,5.00){\line(0,1){10.00}}
\put(25.00,5.00){\line(0,1){10.00}}
\put(40.00,5.00){\line(0,1){10.00}}
\put(45.00,5.00){\line(0,1){10.00}}
\put(60.00,5.00){\line(0,1){10.00}}
\put(80.00,5.00){\line(0,1){10.00}}

\put(10.00,5.00){\line(0,1){5.00}}
\put(15.00,5.00){\line(0,1){5.00}}
\put(30.00,5.00){\line(0,1){5.00}}
\put(35.00,5.00){\line(0,1){5.00}}

\put(50.00,5.00){\line(0,1){5.00}}
\put(55.00,5.00){\line(0,1){5.00}}
\put(65.00,5.00){\line(0,1){5.00}}
\put(70.00,5.00){\line(0,1){5.00}}
\put(75.00,5.00){\line(0,1){5.00}}

\put(10.00,10.00){\line(1,0){5.00}}
\put(30.00,10.00){\line(1,0){5.00}}
\put(50.00,10.00){\line(1,0){5.00}}
\put(65.00,10.00){\line(1,0){10.00}}

\put(-25.00,33.00){\footnotesize ${\rm partition}\;\pi$}
\put(-25.00,8.00){\footnotesize ${\rm refinement}\;\eta$}

\put(12.00,17.00){\footnotesize $\eta_{1}'$}
\put(38.00,17.00){\footnotesize $\eta_{2}'$}
\put(52.00,17.00){\footnotesize $\eta_{3}'$}
\put(78.00,17.00){\footnotesize $\eta_{4}'$}
\put(11.00,2.00){\footnotesize $\pi_{1}''$}
\put(31.00,2.00){\footnotesize $\pi_{2}''$}
\put(51.00,2.00){\footnotesize $\pi_{3}''$}
\put(69.00,2.00){\footnotesize $\pi_{4}''$}

\put(0.00,-15.00){\line(1,0){5.00}}
\put(10.00,-15.00){\line(1,0){5.00}}
\put(20.00,-15.00){\line(1,0){5.00}}
\put(30.00,-15.00){\line(1,0){5.00}}
\put(50.00,-15.00){\line(1,0){5.00}}
\put(65.00,-15.00){\line(1,0){10.00}}
\put(20.00,-15.00){\line(1,0){5.00}}

\put(0.00,-20.00){\line(0,1){5.00}}
\put(5.00,-20.00){\line(0,1){5.00}}
\put(10.00,-20.00){\line(0,1){5.00}}
\put(15.00,-20.00){\line(0,1){5.00}}
\put(20.00,-20.00){\line(0,1){5.00}}
\put(25.00,-20.00){\line(0,1){5.00}}
\put(30.00,-20.00){\line(0,1){5.00}}
\put(35.00,-20.00){\line(0,1){5.00}}
\put(40.00,-20.00){\line(0,1){5.00}}
\put(45.00,-20.00){\line(0,1){5.00}}
\put(50.00,-20.00){\line(0,1){5.00}}
\put(55.00,-20.00){\line(0,1){5.00}}

\put(60.00,-20.00){\line(0,1){5.00}}
\put(65.00,-20.00){\line(0,1){5.00}}
\put(70.00,-20.00){\line(0,1){5.00}}
\put(75.00,-20.00){\line(0,1){5.00}}
\put(80.00,-20.00){\line(0,1){5.00}}

\put(-25.00,-17.00){\footnotesize ${\rm refinement}\;\sigma$}

\put(1.00,-23.00){\footnotesize $\sigma_{1}'$}
\put(21.00,-23.00){\footnotesize $\sigma_{2}'$}
\put(39.00,-23.00){\footnotesize $\sigma_{3}'$}
\put(44.00,-23.00){\footnotesize $\sigma_{4}'$}
\put(59.00,-23.00){\footnotesize $\sigma_{5}'$}
\put(79.00,-23.00){\footnotesize $\sigma_{6}'$}

\put(11.00,-23.00){\footnotesize $\pi_{1}''$}
\put(31.00,-23.00){\footnotesize $\pi_{2}''$}
\put(51.00,-23.00){\footnotesize $\pi_{3}''$}
\put(69.00,-23.00){\footnotesize $\pi_{4}''$}

\end{picture}
\\
\begin{center}
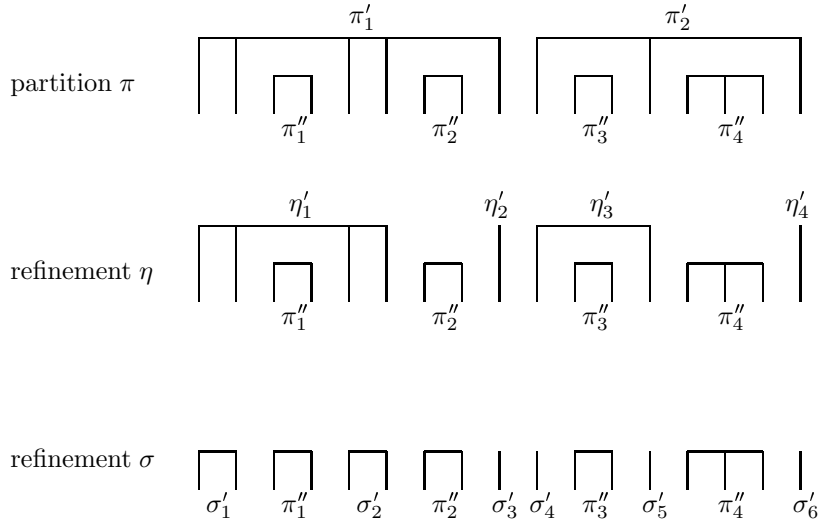
{\it Figure 1.} Example of a decomposition triple $\sigma \preceq \eta \preceq \pi$.\\[15pt]
\end{center}
\indent{\par}
Finally, to given $\eta\in {\cal P}(\pi)$ we associate its coarsest interval subpartition
$\sigma\in {\cal I}(n)$, i.e. $\sigma=\sigma'\cup \sigma''$, where $\sigma''=\eta''$
and $\sigma'$ is the coarsest interval subpartition of $\eta'$.
Note that $\sigma$ is also the coarsest interval subpartition 
of $\pi$ and every block $\pi'_{j}$ gives rise to an odd number of subblocks of $\sigma'$.
In such a way we obtain a {\it decomposition triple}
$$
\sigma\preceq \eta \preceq \pi 
$$
which can be nicely illustrated in terms of diagrams.
We can think of blocks of $\pi'$ as 'bridges' lying above blocks of $\pi''$ (see partition $\pi$ in Figure 1).
Now, $\eta$ is obtained from $\pi$ by erasing certain `bridge connections' over inner blocks
(the latter remain unchanged). Finally, $\sigma$ is obtained from $\eta$ by
erasing the remaining `bridge connections' over inner blocks (the latter remain unchanged).\\
\indent{\par}
{\sc Example 3.2.}
Consider the partition $\pi\in {\cal D}_{2}(17)$ consisting of 2 outer blocks $\pi_{1}'=\{1,2,5,6,9\}$, $\pi_{2}'=\{10,13,17\}$
and 4 inner blocks $\pi_{1}''=\{3,4\}$, $\pi_{2}''=\{7,8\}$, $\pi_{3}''=\{11,12\}$, $\pi_{4}''=\{14,15,16\}$.
Let $\eta$ be its refinement obtained by splitting the block $\pi_{1}'$ into two subblocks: $\eta_{1}'=\{1,2,5,6\}$
and $\eta_{2}'=\{9\}$, and block $\pi_{2}'$ into two subblocks: $\eta_{3}'=\{10,13\}$ and $\eta_{4}'=\{17\}$. Here,
$\eta'=\{\eta_{1}',\eta_{2}',\eta_{3}',\eta_{4}'\}$.
Clearly, $(\pi,\eta)\in {\cal DP}_{2}(17)$ (the pair is shown in Figure 1).
Finally, the coarsest interval subpartition $\sigma$ of $\pi$ is given by
$\sigma=\sigma'\cup \pi''$, where $\sigma'=\{\{1,2\}, \{5,6\},\{9\}, \{10\},\{13\},\{17\}\}$. \\
\indent{\par}
We complete this section with two technical propositions.\\
\indent{\par}
{\sc Proposition 3.1.}
{\it For every $n\in {\mathbb N}$ there is a bijection between ${\cal D}_{2}(n)$ and
the set ${\cal C}(n)$ of pairs $(\tau,\sigma)$, where $\tau\in {\cal I}(n)$ and
$\sigma \in {\cal I}_{{\rm odd}}(\tau)$.}\\
\indent{\par}
{\it Proof.}
Let $\pi\in {\cal D}_{2}(n)$ be given and let $\pi=\pi'\cup \pi''$ be its decomposition (3.6).
Let $f:{\cal D}_{2}(n)\rightarrow {\cal C}(n)$, where
$f(\pi)=(\tau, \sigma)$ is defined as follows.
For every $1\leq j \leq r$, define the block $\tau_{j}$ to be the union of
the block $\pi_{j}$ and all blocks of $\pi$ which are inner w.r.t. $\pi_{j}$.
Then $\sigma$ is defined to be the coarsest interval refinement of
$\pi$ (obtained by `erasing all bridge connections' in $\pi$).
It can be seen that $f$ is a bijection (using diagrams, giving the pair $(\tau, \sigma)$
specifies blocks of $\pi$ in two steps: first we give the intervals which are `covered' by
outer blocks of $\pi$ and then we split up every such interval into an odd number of subintervals
which show the positions of the inner blocks of $\pi$). \hfill
$\blacksquare$\\
\indent{\par}
{\sc Proposition 3.2.}
{\it For every $n\in {\mathbb N}$ there is a bijection between ${\cal DP}_{2}(n)$
and the set ${\cal F}(n)$ of triples $(m,\sigma , j)$, where $1\leq m \leq n$,
$\sigma\in {\cal I}(m)$ and $j=(j_{1},j_{2}, \ldots , j_{m-1})$ is a tuple of non-negative integers
whose sum is equal to $n-m$.}\\
\indent{\par}
{\it Proof.}
The bijection $g: {\cal DP}_{2}(n) \rightarrow {\cal F}(n)$ is given by
$g(\pi,\eta)=(m,\sigma,j)$, where the triple $(m,\sigma,j)$ is defined as follows.
First, we set $m=|\pi_{1}'|+|\pi_{2}'|+\ldots +|\pi_{r}'|$, i.e. $m$ counts all numbers which
belong to the outer blocks of $\pi$. Then we define $\sigma\in {\cal I}(m)$
as the unique partition of the number $m$ which corresponds to the partition $\eta'$ of
the $m$-element set $\pi_{1}'\cup \pi_{2}'\cup \ldots \cup \pi_{r}'$ into blocks of $\eta'$.
Finally, we set $j=(j_{1},j_{2}, \ldots , j_{m-1})$, where $j_{k}$ is the nonnegative integer
equal to the size of the inner block of $\pi$ which immediately follows the $k$-th leg of $\pi'$
($\pi'$ has $m$ legs but the last leg of $\pi'$ ends the diagram, so it does not count).
Of course, if there is no inner block following the $k$-th leg of $\pi'$,
then we set $j_{k}=0$. It can be seen that the mapping $g$ is a bijection.
Using diagrams, one can say that by giving the triple $(m,\sigma,j)$, we simply draw
the diagram corresponding to $(\pi,\eta)$ in the following order: first we draw all
outer blocks of $\eta'$ and then every inner block of $\eta''=\pi''$
is drawn on the right side of the suitable leg of $\eta'$.\hfill $\blacksquare$\\[15pt]
\myownsection
\begin{center}
{\sc 4. Orthogonal structures}
\end{center}
In this Section we introduce the notion of `orthogonal subalgebras' of a given (*- , $C^{*}$-) algebra with respect
to a pair of functionals (states) and the corresponding notions of the `orthogonal product' of two Hilbert
spaces and the `orthogonal product' of two (*-, $C^{*}$-) algebras. We then construct `orthogonal random variables'
with prescribed probability distributions. To some extent, these structures resemble the corresponding monotone structures and for that reason can be viewed as `quasi-monotone'.

The orthogonal product is neither commutative nor associative, but it turns out
useful in the construction of decompositions of
the free additive convolution of measures.\\
\indent{\par}
{\sc Definition 4.1.}
Let $({\cal A},\varphi, \psi)$ be a unital algebra with a pair
of linear normalized functionals and let ${\cal A}_{1}$ and ${\cal A}_{2}$
be non-unital subalgebras of ${\cal A}$.
We say that ${\cal A}_{2}$ is {\it orthogonal} to ${\cal A}_{1}$
with respect to $(\varphi, \psi)$ if\\
\indent{\par}
(i)
$\;\varphi(ba_{2})=\varphi(a_{1}b)=0$
\indent{\par}
(ii)
$\varphi(w_{1}a_{1}ba_{2}w_{2})=\psi(b)
\left(\varphi(w_{1}a_{1}a_{2}w_{2})- \varphi(w_{1}a_{1})\varphi(a_{2}w_{2})\right)$\\[5pt]
for any $a_{1},a_{2}\in {\cal A}_{1}$, $b\in {\cal A}_{2}$
and any elements $w,v$ of the algebra ${\rm alg}({\cal A}_{1},{\cal A}_{2})$
generated by ${\cal A}_{1}$ and ${\cal A}_{2}$.
We say that the pair $(a,b)$ of elements of ${\cal A}$ is {\it orthogonal} with respect to $(\varphi, \psi)$
if the algebra generated by $a\in {\cal A}$ is orthogonal to the algebra generated by
$b\in {\cal A}$ .\\
\indent{\par}
{\sc Remark 4.1.}
Note that $\varphi$ is uniquely determined on the algebra generated by ${\cal A}_{1}$ and ${\cal A}_{2}$
by restrictions $\varphi|{\cal A}_{1}$ and $\psi|{\cal A}_{2}$.
In the case of $\psi$, the situation is quite different. In fact,
this is only for the sake of convenience that we consider two states
$\varphi, \psi$ on all of ${\cal A}$ (it is natural in the Hilbert
space setting, where we choose states associated with unit vectors of the `large' Hilbert space).
For our purposes, it would be sufficient to assume $\psi$ to be defined 
only on the subalgebra ${\cal A}_{2}$. Another observation is that `orthogonality' w.r.t.
$(\varphi, \psi)$, is quite different from `conditional freeness' w.r.t. $(\varphi, \psi)$
studied in [7], although it also involves two states on ${\cal A}$. \\
\indent{\par}
Let us begin with the Hilbert space setting and introduce the notion of
an orthogonal product of two Hilbert spaces with distinguished unit vectors.\\
\indent{\par}
{\sc Definition 4.2.}
Let $({\cal H}_{1}, \xi_{1})$ and $({\cal H}_{2}, \xi_{2})$ be
Hilbert spaces with distinguished unit vectors $\xi_{1}$ and $\xi_{2}$, respectively.
The {\it orthogonal product} of $({\cal H}_{1}, \xi_{1})$ and $({\cal H}_{2}, \xi_{2})$
is the pair $({\cal H}, \xi)$, where
\begin{equation}
{\cal H}={\mathbb C}\xi \oplus {\cal H}_{1}^{0}\oplus ({\cal H}_{2}^{0}\otimes {\cal H}_{1}^{0}),
\end{equation}
with ${\cal H}_{i}^{0}={\cal H}_{i}\ominus {\mathbb C}\xi_{i}$ denoting  the orthogonal complement
of ${\mathbb C}\xi_{i}$, $i=1,2$ and $\xi$ being a unit vector.
We denote it by $({\cal H}, \xi)=({\cal H}_{1}, \xi_{1})\vdash ({\cal H}_{2}, \xi_{2})$
and by $\varphi$ - the canonical state on ${\cal B}({\cal H})$ associated with the vector $\xi$.\\
\indent{\par}
Note that the orthogonal product of Hilbert spaces is slightly smaller than their monotone
product [20]. In fact, the monotone product of $({\cal H}_{1}, \xi_{1})$ and
$({\cal H}_{2}, \xi_{2})$ is equal to the direct sum of their orthogonal product
and ${\cal H}_{2}^{0}$.
Clearly, $({\cal H}_{1}, \xi_{1})\vdash ({\cal H}_{2}, \xi_{2})$ is also 
a truncation of the free product
of Hilbert spaces $({\cal H}_{1}, \xi_{1})*({\cal H}_{2}, \xi_{2})$ [26].
However, in order to study representations,
it is more convenient to use the tensor product
${\cal H}_{1}\otimes {\cal H}_{2}$ and an isometry
$U: {\cal H}\rightarrow {\cal H}_{1}\otimes {\cal H}_{2}$ given by
\begin{equation}
U(\xi)=\xi_{1}\otimes \xi_{2},\;\; U(h_{1})=h_{1}\otimes \xi_{2},\;\; U(h_{2}\otimes h_{1})=h_{1}\otimes h_{2}
\end{equation}
for any $h_{1}\in {\cal H}_{1}^{0}$ and $h_{2}\in {\cal H}_{2}^{0}$. In particular, we have
\begin{equation}
UU^{*}=1-P_{{\mathbb C}\xi_{1}\otimes {\cal H}_{2}^{0}}.
\end{equation}
Using the isometry $U$, we shall define *-representations
$\tau_{i}:{\cal B}({\cal H}_{i})\rightarrow {\cal B}({\cal H})$
by
\begin{equation}
\tau_{1}(a)= U^{*}(a\otimes P_{2}) U, \;\;\;
\tau_{2}(b)= U^{*}(P_{1}^{\perp} \otimes b)U
\end{equation}
where $P_{1}$, $P_{2}$ are the projections
onto ${\mathbb C}\xi_{1}$ and ${\mathbb C}\xi_{2}$, respectively.
Note that $\tau_{1}$ and $\tau_{2}$ are faithful non-unital *-homomorphisms.

In the theorem below we describe the properties of the
orthogonal product of Hilbert spaces.\\
\indent{\par}
{\sc Theorem 4.1.}
{\it Let $({\cal H}, \xi)= ({\cal H}_{1},\xi_{1})\vdash ({\cal H}_{2}, \xi_{2})$ be the
orthogonal product of Hilbert spaces and let $\varphi$, $\varphi_{1}$ and $\varphi_{2}$ be
the states associated with $\xi$, $\xi_{1}$ and $\xi_{2}$, respectively.
Moreover, let $\psi$ be the state on ${\cal B}({\cal H})$
associated with any unit vector $\eta\in {\cal H}_{1}^{0}\subset {\cal H}$.
Then
\indent{\par}
(i) \;\;${\cal A}_{2}=\tau_{2}({\cal B}({\cal H}_{2}))$
is orthogonal to ${\cal A}_{1}=\tau_{1}({\cal B}({\cal H}_{1}))$ with respect to
$(\varphi, \psi)$,
\indent{\par}
(ii) \;$\varphi \circ \tau_{1}$ agrees with the expectation $\varphi_{1}$ on ${\cal B}({\cal H}_{1})$
\indent{\par}
(iii) $\psi \circ \tau_{2}$ agrees with the expectation $\varphi_{2}$ on ${\cal B}({\cal H}_{2})$.}\\
\indent{\par}
{\it Proof.}
Let $a\in {\cal A}_{1}$, $b\in {\cal A}_{2}$ and $w_{1}\in {\rm alg}({\cal A}_{1},{\cal A}_{2})$.
First, observe that (4.3) implies that
$$
\varphi(w_{1}\tau_{2}(b))=\varphi(w_{1}U^{*}(P_{1}^{\perp}\otimes b)U)=0
$$
since $(P_{1}^{\perp}\otimes b)(\xi_{1}\otimes \xi_{2})=0$.
Similarly, $\varphi(\tau_{2}(b)w_{1})=0$ and thus the
condition (i) of Definition 4.1 holds. We need to show the condition (ii) of that
definition.
We have
\begin{eqnarray*}
\varphi(w_{1}\tau_{1}(a_{1})\tau_{2}(b)\tau_{1}(a_{2})w_{2}) &=&
\varphi\left(
w_{1}U^{*}(a_{1}\otimes P_{2})UU^{*}(P_{1}^{\perp}\otimes b)UU^{*}(a_{2}\otimes P_{2})Uw_{2}
\right)\\
&=&
\varphi\left(
w_{1}U^{*}(a_{1}\otimes P_{2})UU^{*}(P_{1}^{\perp}a_{2}\otimes bP_{2})Uw_{2}
\right)
\end{eqnarray*}
for any $w_{1},w_{2}\in {\rm alg}({\cal A}_{1},{\cal A}_{2})$
in view of (4.3). For the same reason, we get
$$
(a_{1}\otimes P_{2})UU^{*}(P_{1}^{\perp}a_{2}\otimes bP_{2})Uw_{2}=\varphi_{2}(b)(a_{1}\otimes P_{2})
UU^{*}(P_{1}^{\perp}a_{2}\otimes P_{2})Uw_{2}.
$$
Below we shall demonstrate that $\varphi_{2}(b)=\psi(\tau_{2}(b))$ for the state
$\psi$ associated with any unit vector
$\eta$ from ${\cal H}_{1}^{0}$. Therefore, we are left with computing
$$
\varphi\left(
w_{1}U^{*}(a_{1}\otimes P_{2})UU^{*}(P_{1}^{\perp}a_{2}\otimes P_{2})Uw_{2}
\right)
$$
\begin{eqnarray*}
&=&
\varphi\left(
w_{1}U^{*}(a_{1}\otimes P_{2})UU^{*}(a_{2}\otimes P_{2})Uw_{2}
\right)
-
\varphi\left(
w_{1}U^{*}(a_{1}\otimes P_{2})UU^{*}(P_{1}a_{2}\otimes P_{2})Uw_{2}
\right)\\
&=&
\varphi\left(
w_{1}U^{*}(a_{1}\otimes P_{2})UU^{*}(a_{2}\otimes P_{2})Uw_{2}
\right)
-
\varphi\left(
w_{1}U^{*}(a_{1}\otimes P_{2})UP_{\xi}U^{*}(a_{2}\otimes P_{2})Uw_{2}
\right)\\
&=&
\varphi(w_{1}a_{1}a_{2}w_{2})-
\varphi(w_{1}a_{1})\varphi(a_{2}w_{2})
\end{eqnarray*}
where we used $U^{*}(P_{1}\otimes P_{2})= P_{\xi}U^{*}(1\otimes P_{2})$, with
$P_{\xi}$ denoting the projection onto ${\mathbb C}\xi$.
Finally,
\begin{eqnarray*}
\varphi\circ \tau_{1}(a)&=&\langle (a\otimes P_{2})
\xi_{1}\otimes \xi_{2}, \xi_{1}\otimes \xi_{2}\rangle = \varphi_{1}(a)\\
\psi\circ \tau_{2}(b)&=&\langle  (P_{1}^{\perp}\otimes b)\eta \otimes \xi_{2}, \eta \otimes \xi_{2}\rangle = \varphi_{2}(b)
\end{eqnarray*}
for any $a\in {\cal B}({\cal H}_{1})$ and $b\in {\cal B}({\cal H}_{2})$, which completes the proof.
\hfill $\blacksquare$\\
\indent{\par}
{\sc Corollary 4.2.}
{\it Let $\mu, \nu$ be compactly supported probability measures on ${\mathbb R}$.
Then there exist a Hilbert space ${\cal H}$, unit vectors $\xi, \eta\in {\cal H}$ and
self-adjoint bounded random variables $X_1,X_2\in {\cal B}({\cal H})$ such that the pair $(X_1,X_2)$ is orthogonal
w.r.t. $(\varphi, \psi)$, where $\varphi$ and $\psi$ are vector states associated with
$\xi, \eta\in {\cal H}$. Moroever, the $\varphi$-distribution of $X_1$ and the
$\psi$-distribution of $X_2$ coincide with $\mu$ and $\nu$, respectively. Finally,
the $\varphi$-distribution of $X_1+X_2$, denoted $\mu \vdash \nu$, is compactly supported.}\\
\indent{\par}
{\it Proof.}
Let ${\cal H}_{1}=L^{2}({\mathbb R}, \mu)$ and ${\cal H}_{2}=L^{2}({\mathbb R}, \nu)$
and take $\xi_{1}=1$ and $\xi_{2}=1$. Let $\hat{x}_{1}$ and $\hat{x}_{2}$ be the
standard multiplication operators on these spaces, namely
$\hat{x}_{1}f(x_1)=x_1f(x_1)$ and $\hat{x}_{2}g(x_2)=x_2g(x_2)$. They are bounded self-adjoint operators with
distributions $\mu$ and $\nu$, respectively. By taking the orthogonal product
$({\cal H}, \xi)= ({\cal H}_{1}, \xi_{1})\vdash ({\cal H}_{2}, \xi_{2})$
we can construct bounded self-adjoint random variables $X_1=\tau_{1}(\hat{x}_1)$ and $X_2=\tau_{2}(\hat{x}_2)$
from ${\cal B}({\cal H})$
such that the pair $(X_1,X_2)$ is orthogonal w.r.t. $(\varphi, \psi)$, where $\varphi$ is the
vector state on ${\cal B}({\cal H})$ associated with $\xi$, and $\psi$ is the
vector state on ${\cal B}({\cal H})$ associated with any function
$f\in L^{2}({\mathbb R}, \mu)$ which satisfies $\int_{\mathbb R}f(x_1)\mu(dx_1)=0$
and $\int_{\mathbb R}f^{2}(x_1)\mu(dx_1)=1$.
Finally, it is clear that the sum $X_1+X_2$ is a bounded self-adjoint operator on ${\cal H}$ and thus
its probability distribution extends to a compactly supported measure on the real line.
\hfill $\blacksquare$\\
\indent{\par}
{\sc Example 4.1.}
Let $({\cal G}_{1},e_{1})$ and $({\cal G}_{2},e_{2})$ be two uniformly locally 
finite rooted graphs with adjacency matrices $A_{1}$ and $A_{2}$, respectively, which extend to
bounded operators on ${\cal H}_{1}=l_{2}(V_{1})$ and ${\cal H}_{2}=l_{2}(V_{2})$, 
where $V_{1}$ and $V_{2}$ denote their sets of vertices. Then 
$$
A^{(1)}=A_{1}\otimes P_{e_{2}}\,\,\,\,{\rm and}\,\,\,\,A^{(2)}=P_{e_{1}}^{\perp}\otimes A_{2}
$$
are orthogonal w.r.t. $(\varphi, \psi)$, where $\varphi$ and $\psi$ are states
on ${\cal B}({\cal H}_{1})\otimes {\cal B}({\cal H}_{2})$ associated with
vectors $\delta(e_{1})\otimes \delta (e_{2})$ and $\delta (v)\otimes \delta (e_{2})$,
respectively, with $v\in V_{1}^{0}=V_{1}\setminus \{e_{1}\}$, and $P_{e_{i}}$ 
is the projection onto ${\mathbb C}\delta (e_{i})$. The sum $A=A^{(1)}+A^{(2)}$
is the adjacency matrix of a uniformly locally finite rooted graph 
$({\cal G}_{1}\vdash {\cal G}_{2},e)$
obtained by attaching a replica of ${\cal G}_{2}$ by its root to every vertex of
$V_{1}^{0}$ and setting $e=e_{1}\times e_{2}$, called
the {\it orthogonal product} of $({\cal G}_{1},e_{1})$ and $({\cal G}_{2},e_{2})$.
The matrix $A$ extends to a bounded operator on the orthogonal product
$({\cal H}_{1}, \delta(e_{1}))\vdash ({\cal H}_{2}, \delta(e_{2}))$.
By Corollary 4.2, the spectral distribution of $A$ in the state associated with vector $\delta(e)$
is given by $\mu \vdash \nu$, where $\mu$ and $\nu$ are spectral distributions of
$A_{1}$ and $A_{2}$ associated with $\delta(e_{1})$ and $\delta(e_{2})$, 
respectively. 
A detailed study of the orthogonal product of rooted graphs will be given in a separate paper
[1].
\\
\indent{\par}
It is now natural to define the orthogonal product in the setting of $C^{*}$-probability spaces.
If $({\cal A}_{i}, \varphi_{i})$, $i=1,2$, are $C^{*}$-probability spaces and
$({\cal H}_{i}, \pi_{i}, \xi_{i})$ - the corresponding GNS triples, we first construct
the orthogonal product of Hilbert spaces
$({\cal H}, \xi)=({\cal H}_{1}, \xi_{1})\vdash ({\cal H}_{2}, \xi_{2})$
and then define representations
$$
\iota_{i}: {\cal A}_{i}\rightarrow {\cal B}({\cal H}), \;\;\;
\iota_{i}=\tau_{i}\circ \pi_{i}
$$
where $\tau_{i}$, $i=1,2$, are given by (4.4).
Let ${\cal A}$ be the $C^{*}$-algebra generated by subalgebras $\iota_{1}({\cal A}_{1})$
and $\iota_{2}({\cal A}_{2})$ of ${\cal B}({\cal H})$ and the identity $I\in {\cal B}({\cal H})$ and let
$\varphi$ denote the state on ${\cal A}$ associated with
the vector $\xi$. Then the pair $({\cal A}, \varphi)$ is
called the {\it orthogonal product} of $C^{*}$-probability spaces
$({\cal A}_{1}, \varphi_{1})$ and $({\cal A}_{2}, \varphi_{2})$ and is denoted
$({\cal A}_{1}, \varphi_{1})\vdash ({\cal A}_{2}, \varphi_{2})$.
\\
\indent{\par}
{\sc Theorem 4.3.}
{\it Let $({\cal A}, \varphi)=({\cal A}_{1}, \varphi_{1})\vdash ({\cal A}_{2}, \varphi_{2})$
be the orthogonal product of $C^{*}$-probability spaces equipped with the natural *-homomorphisms
$\iota_{i}:{\cal A}_{i}\rightarrow {\cal A}$ and let $\psi$ be the state on ${\cal A}$
associated with any unit vector $\eta\in {\cal H}_{1}^{0}$. Then
\indent{\par}
(i) \; $\iota_{1}({\cal A}_{1})$ is orthogonal to $\iota_{2}({\cal A}_{2})$
w.r.t. $(\varphi, \psi)$,
\indent{\par}
(ii) \;$\varphi \circ \iota_{1}$  agrees with the expectation $\varphi_{1}$ on ${\cal A}_{1}$
\indent{\par}
(iii) $\psi \circ \iota_{2}$ agrees with the expectation $\varphi_{2}$ on ${\cal A}_{2}$.}\\
\indent{\par}
{\it Proof.}
This is Theorem 4.1 adapted to the $C^{*}$-algebra setting (for an analogous formulation
in the monotone case, see [20]). \hfill $\blacksquare$\\
\indent{\par}
{\sc Example 4.2.}
In Example 4.1, let ${\cal A}_{i}$ be the $C^{*}$-algebra generated by $A_{i}$ 
and the identity $I_{i}$ on ${\cal H}_{i}$ and let $\varphi_{i}$ be the state on ${\cal A}_{i}$ asssociated
with the vector $\delta (e_{i})$, $i=1,2$.
Then the pair $({\cal A}, \varphi)$, where ${\cal A}$ is the $C^{*}$-algebra generated by 
$A^{(1)}$ and $A^{(2)}$ and the identity $I_{1}\otimes I_{2}$ and $\varphi$ is the state 
on ${\cal A}\subset {\cal B}({\cal H}_{1})\otimes {\cal B}({\cal H}_{2})$ associated
with the vector $\delta(e_{1})\otimes \delta(e_{2})$, is the orthogonal product of $({\cal A}_{1},\varphi_{1})$ and $({\cal A}_{2}, \varphi_{2})$.\\
\indent{\par}
{\sc Remark 4.2.}
The notion of the orthogonal product can also be introduced in the category
of noncommutative (*-) probability spaces. Then, conditions of
Definition 4.1 can be used as defining conditions for 
the orthogonal product of functionals (states)
on the free product ${\cal A}_{1}\sqcup{\cal A}_{2}$
without identification of units. One can use extensions $\widetilde{\cal A}_{1}={\cal A}_{1}*{\mathbb C}[p_{1}]$
and $\widetilde{\cal A}_{2}={\cal A}_{2}*{\mathbb C}[p_{2}]$ by idempotents (projections) $p_{1}$ and $p_{2}$ 
to construct the unital (*-) homomorphism
$j: {\cal A}_{1}\sqcup {\cal A}_{2}\rightarrow \widetilde{\cal A}_{1}\otimes \widetilde{\cal A}_{2}$
as the linear and multiplicative extension of
$j(a)=a\otimes p_{2}$ and $j(b)=p_{1}^{\perp}\otimes b$
for any $a\in {\cal A}_{1}$ and $b\in {\cal A}_{2}$, where $p_{1}^{\perp}=1-p$.
Then $\varphi$ agrees with the functional (state) 
$(\widetilde{\varphi}_{1}\otimes  \widetilde{\phi}_{2})\circ j$.
In particular, in the case of *-probability spaces, this
proves positivity of $\varphi$. Therefore, the pair $({\cal A}_{1}\sqcup {\cal A}_{2},
\varphi)$ can be defined as the orthogonal product of noncommutative (*)-probability spaces
$({\cal A}_{1}, \varphi_{1})$ and $({\cal A}_{2}, \varphi_{2})$.
\\[10pt]
\myownsection
\begin{center}
{\sc 5. Orthogonal convolution}
\end{center}
The moments of the `orthogonal convolution' $\mu \vdash \nu$ of compactly supported 
probability measures can be computed using the Hilbert space realization of Section 4.
Keeping the notations of Corollary 4.2, we obtain the following proposition. \\
\indent{\par}
{\sc Proposition 5.1.}
{\it For any $\pi\in {\cal I}(n)$ it holds that}
\begin{equation}
m_{\mu \,\vdash \,\nu}(\pi)=\sum _{\sigma \in {\cal I}_{{\rm odd}}(\pi)}
(-1)^{|\sigma'|-|\pi|}k_{\mu}^{*}(\sigma ')m_{\nu} (\sigma '')
\end{equation}
{\it where $k_{\mu}^{*}(\sigma')$ and $m_{\nu}(\sigma'')$ are given by (3.3) and (3.5) and
$\sigma=\sigma'\cup \sigma''$ is the decomposition given by (3.1).}\\
\indent{\par}
{\it Proof.}
First consider the case, when $\pi$ consists of
one $n$-element block, which we denote $\pi=(n)$, where $n\in {\mathbb N}$. We have
\begin{eqnarray*}
(\mu \vdash \nu)(n)
&=&
\langle (X_1\otimes P_{2} + P_{1}^{\perp}\otimes X_2)^{n} \xi_{1}\otimes \xi_{2},
\xi_{1}\otimes \xi_{2}\rangle\\
&=&
\sum_{r=1}^{n}
\sum_{
\stackrel{j_{1}+j_{2}+ \ldots + j_{r}=n}
{\scriptscriptstyle r\;{\rm odd}}
}
\langle
X_1^{j_{1}}P_{1}^{\perp}X_1^{j_{3}}P_{1}^{\perp}\ldots P_{1}^{\perp}X_1^{j_{r}}\xi_{1},\xi_{1}
\rangle\\
&&
\;\;\;\;\;\;\;\;\;\;\;\;\;\;\;\;\;\times \;\;
\langle
P_{2}X_2^{j_{2}}P_{2}\ldots X_2^{j_{r-1}}P_{2} \xi_{2}, \xi_{2}
\rangle
\end{eqnarray*}
where we understand that the summation runs over the set of ordered partitions of the number $n$
and thus all $j_{k}$'s are assumed to be non-zero.
Now, if we denote by $\sigma$ the interval partition associated with
the tuple $(j_{1},j_{2}, \ldots , j_{n})$, and by $\sigma=\sigma'\cup \sigma''$ - the alternating
decomposition (3.1), we get
$$
\langle
X_1^{j_{1}}P_{1}^{\perp}X_1^{j_{3}}P_{1}^{\perp}\ldots P_{1}^{\perp}X_1^{j_{r}}\xi_{1},\xi_{1}
\rangle
=
\sum_{\eta \succeq \sigma'} (-1)^{|\eta|-1}m_{\mu}(\eta)= (-1)^{|\sigma'|-1}k_{\mu}^{*}(\sigma ')
$$
where we use (3.5). Moroever,
$$
\langle
P_{2}X_2^{j_{2}}P_{2}\ldots X_2^{j_{r-1}}P_{2} \xi_{2}, \xi_{2}
\rangle
=
\nu(j_{1})\nu(j_{2})\ldots \nu(j_{r-1})\\
=
m_{\nu}(\sigma'').
$$
This gives (5.1) for $\pi=(n)$. It remains to extend this result multiplicatively to any $\pi\in {\cal I}(n)$.
Namely, (5.1) holds for every block $\pi_{j}$ of $\pi$ with summation running over partitions
$\sigma(j)\in {\cal I}_{{\rm odd}}(|\pi_{j}|)$ with the sign factor equal to $(-1)^{|\sigma'(j)-1|}$.
Thus, every block $\pi_{j}$ is decomposed into and odd number of subblocks. The sum over $\sigma(j)$'s
gives (5.1) for any $\pi$ since $\sum_{j}(|\sigma'(j)|-1)=|\sigma'|-|\pi|$.
\hfill $\blacksquare$\\
\indent{\par}
One can generalize the definition of the orthogonal convolution and Proposition 5.1
to distributions of an arbitrary orthogonal pair $(a,b)$ of elements of 
a noncommutative probability space ${\cal A}$ (see Remark 4.2). \\
\indent{\par}
{\sc Definition 5.1.}
Let $(a,b)$ be a pair of random variables from a unital algebra ${\cal A}$
which is orthogonal w.r.t. to a pair of normalized linear functionals $(\varphi, \psi)$,
with $\mu$ denoting the $\varphi$-distribution of $a$ and $\nu$ denoting the $\psi$-distribution
of $b$. By the {\it orthogonal convolution} $\mu \vdash \nu$ we understand the
$\varphi$-distribution of $a+b$.\\
\indent{\par}
{\sc Example 5.1.}
Using Proposition 5.1 as well as (3.3) and (3.5), we get
\begin{eqnarray*}
\mu_{a+b}(1)&=&\mu_{a}(1)\\
\mu_{a+b}(2)&=&\mu_{a}(2)\\
\mu_{a+b}(3)&=&\mu_{a}(3)+ (\mu_{a}(2)-\mu_{a}^{2}(1))\nu_{b}(1)\\
\mu_{a+b}(4)&=& \mu_{a}(4)+2\mu_{a}(3)\nu_{b}(1)+\mu_{a}(2)\nu_{b}(2)\\
&&-2\mu_{a}(2)\mu_{a}(1)\nu_{b}(1)-\mu_{a}^{2}(1)\nu_{b}(2)
\end{eqnarray*}
It can be seen that $\mu_{a+b}$ is not symmetric with respect to $\mu_{a}$ and $\nu_{b}$.
In particular, the first two moments of $a+b$ agree with the moments of
$a$. \\
\indent{\par}
More generally, the moment $\mu_{a+b}(n)$ for $n\geq 2$ can be expressed in terms
of the moments $\mu_{a}(k)$ of orders $k\leq n$ and the moments $\mu_{b}(l)$ of orders $l\leq n-2$.
In fact, using the language of `universal polynomials',
we obtain the following analogue of Proposition 4.3 of [1] or Proposition 1.2 of [24].\\
\indent{\par}
{\sc Proposition 5.2.}
{\it Let $(a,b)$ be a pair of elements of a unital algebra ${\cal A}$ which
is orthogonal w.r.t. $(\varphi, \psi)$.
The $\varphi$-distribution of $a+b$ depends only on $\mu_{a}$ and $\nu_{b}$
and there are universal polynomials with integer coefficients
$P_{m}(x_{1},\ldots , x_{m}, y_{1}, \ldots , y_{m-2})$ for $m>2$, with
$P_{1}(x_{1})=x_{1}$ and $P_{2}(x_{1},x_{2})=x_{2}$, such that
\indent{\par} (1) $P_{m}$ is homogenous of degree $m$ in the $x$ and $y$ variables taken together,
where degree $j$ is assigned to $x_{j}$ and $y_{j}$,
\indent{\par} (2) $\mu_{a+b}(m)=P_{n}(\mu_{a}(1), \ldots, \mu_{a}(m), \nu_{b}(1), \ldots , \nu_{b}(m-2))$.}\\
\indent{\par}
{\it Proof.}
It follows directly from (5.1) for the partition $\pi=(m)$ consisting of one block
that $\mu_{a+b}$ depends only on $\mu_{a}$ and $\nu_{b}$
since $k_{\mu}^{*}(\sigma')$ and $m_{\nu}(\sigma'')$ depend only on $\mu_{a}$ and $\nu_{b}$.
Now, each $k_{\mu}^{*}(\sigma')$ is a polynomial in the
moments  of $\mu$  with integer coefficients and $m_{\nu}$ is just a product of moments of $\nu$
as the proof of Proposition 5.1 demonstrates. Assigning the variable $x_{j}$ to $\mu(j)$
and $y_{j}$ to $\nu(j)$, we obtain the polynomial $P_{m}$. Since in the expression on the right-hand side
of (5.1) we have a summation over odd subpartitions of $\pi=(m)$, moments $\mu(j)$ of orders $j\leq m$
and $\nu(j)$ of orders $j\leq m-2$ appear and that is why $P_{m}$ depends on $x_{1}, \ldots , x_{m}$
and $y_{1}, \ldots , y_{m-2}$.\hfill $\blacksquare$\\
\indent{\par}
{\sc Corollary 5.3.}
{\it Let $a_{1},a_{2}, \ldots , a_{n}\in {\cal A}$ and let $\mu_{1}, \mu_{2}, \ldots , \mu_{n}$
be their distributions w.r.t. normalized linear functionals $\varphi_{1}, \varphi_{2}, \ldots, \varphi_{n}$
on ${\cal A}$, respectively. If the pair $(a_{j}, a_{j+1}+\ldots + a_{n})$ is orthogonal
w.r.t. $(\varphi_{j}, \varphi_{j+1})$ for every $1\leq j \leq n-1$, then the $m$-th moment of
$a_{1}+a_{2}+\ldots +a_{n}$ depends only on $\mu_{r}(j_{r})$, where $1\leq j_{r} \leq m-2r+2$
and $1\leq r\leq n$.}\\
\indent{\par}
{\it Proof.}
A repeated application of Proposition 5.2 gives the assertion.
\hfill $\blacksquare$\\[10pt]
\myownsection
\begin{center}
{\sc 6. Reciprocal Cauchy transforms}
\end{center}
Our goal now is to derive a formula expressing the reciprocal Cauchy transform of distribution
$\mu\vdash \nu$ in terms of those of $\mu$ and $\nu$.
In the theorem below we shall do it on the level of formal power series for measures
with finite moments of all orders. However, we will also show that
the RHS of (6.2) is the reciprocal Cauchy transform of
a probability measure if $\mu$ and $\nu$ are arbitrary probability measures, which gives
an analytic approach to the orthogonal convolution.

Let us first prove an elementary
combinatorial formula for the reciprocal Cauchy transform $F_{\mu}(z)$
of distribution $\mu$.\\
\indent{\par}
{\sc Proposition 6.1.}
{\it Let $\mu$ and $\nu$ be probability measures with finite moments of all orders.
Then the reciprocal Cauchy transform of the distribution $\mu$ satisfies the
equation}
\begin{equation}
F_{\mu}(z)-z =
\sum_{n=1}^{\infty}\sum_{\pi\in {\cal I}(n)}(-1)^{|\pi|}m_{\mu}(\pi)z^{-n+1}
\end{equation}
{\it where the right-hand side is understood as a formal power series, where
$m_{\mu}(\pi)$ is given by (3.2)-(3.3).}\\
\indent{\par}
{\it Proof.}
We have
\begin{eqnarray*}
F_{\mu}(z)
&=&
\frac{z}{1+\sum_{n=1}^{\infty}\mu(n)z^{-n}}\\
&=&
z(1+\sum_{k=1}^{\infty}(-\sum_{n=1}^{\infty}\mu(n)z^{-n})^{k})\\
&=&
z(1+\sum_{m=1}^{\infty}
\sum_{r=1}^{m}\sum_{j_{1}+j_{2}+\ldots +j_{r}=m}(-1)^{r}\mu(j_{1})\mu(j_{2})\ldots \mu(j_{r})z^{-m})\\
&=&
z+ \sum_{m=1}^{\infty}\sum_{\pi\in {\cal I}(m)}(-1)^{|\pi|}m_{\mu}(\pi)z^{-m+1}
\end{eqnarray*}
which completes the proof.
\hfill $\blacksquare$\\[10pt]
\indent{\par}
{\sc Theorem 6.2.}
{\it Let $\mu$ and $\nu$ be probability measures with finite moments of all orders.
The reciprocal Cauchy transform of $\mu \vdash \nu$ is given by the formula}
\begin{equation}
F_{\mu\,\vdash\, \nu}(z) =F_{\mu}(F_{\nu}(z))-F_{\nu}(z) +z
\end{equation}
{\it where the right-hand side is understood as a formal power series.}\\
\indent{\par}
{\it Proof.}
Using Proposition 6.1, we obtain
\begin{eqnarray*}
L:&=&F_{\mu\,\vdash\,\nu}(z)-z
=
\sum_{n=1}^{\infty}\sum_{\pi\in {\cal I}(m)}(-1)^{|\pi|}m_{\mu \,\vdash \,\nu}(\pi)z^{-n+1}\\
R:&=&F_{\mu}(F_{\nu}(z))-F_{\nu}(z)
=
\sum_{m=1}^{\infty}\sum_{\pi\in {\cal I}(m)}(-1)^{|\pi|}m_{\mu}(\pi)(G_{\nu}(z))^{m-1}.
\end{eqnarray*}
In turn, the definition of $G_{\nu}(z)$ gives
\begin{eqnarray*}
(G_{\nu}(z))^{m-1}&=&
\sum_{j_1=0}^{\infty}\sum_{j_2=0}^{\infty}\ldots \sum_{j_{m-1}=0}^{\infty}
\nu(j_1)\nu(j_2)\ldots \nu(j_{m-1})z^{-j_1-j_2-\ldots -j_{m-1}-(m-1)}\\
&=&
z^{-m+1}
\sum_{k=0}^{\infty}
\left(
\sum_{
\stackrel{j_{1}+j_{2}+\ldots +j_{m-1}=k}
{\scriptscriptstyle j_{1}, j_{2}, \ldots , j_{m-1}\geq 0}
}
\nu(j_1)\nu(j_{2})\ldots \nu(j_{m-1})
\right)
z^{-k}
\end{eqnarray*}
Writing $L$ and $R$ in the form of formal power series
$$
L=\sum_{n=1}^{\infty}L_{n}z^{-n+1}, \;\;\;\;
R=\sum_{n=1}^{\infty}R_{n}z^{-n+1}
$$
we get
\begin{eqnarray*}
L_{n}&=&
\sum_{\tau\in {\cal I}(n)}(-1)^{|\tau|}m_{\mu\,\vdash \,\nu}(\tau)\\
R_{n}&=&
\sum_{m=1}^{n}
\sum_{\sigma\in {\cal I}(m)}(-1)^{|\sigma|}m_{\mu}(\sigma)
\sum_{
\stackrel{j_{1}+j_{2}+\ldots +j_{m-1}=n-m}
{\scriptscriptstyle j_{1}, j_{2}, \ldots , j_{m-1}\geq 0}
}
\nu(j_1)\nu(j_{2})\ldots \nu(j_{m-1})
\end{eqnarray*}
and we thus need to show that $L_{n}=R_{n}$ for every $n\geq 1$. Using Proposition 5.1, formula (3.5)
and Proposition 3.1, we get
\begin{eqnarray*}
L_{n}&=&\sum_{\tau\in {\cal I}(n)}(-1)^{|\tau|}
\sum_{\sigma \in {\cal I}_{{\rm odd}}(\tau)}
(-1)^{|\sigma'|-|\tau|}
k_{\mu}^{*}(\sigma') m_{\nu}(\sigma '')\\
&=&
\sum_{\tau\in {\cal I}(n)}
\sum_{\sigma \in {\cal I}_{{\rm odd}}(\tau)}
\sum_{\eta' \succeq \sigma'}
(-1)^{|\eta'|}m_{\mu}(\eta') m_{\nu}(\sigma '')\\
&=&
\sum_{\pi\in {\cal D}_{2}(n)}\sum_{\eta' \succeq \pi'}
m_{\mu}(\eta')m_{\nu}(\pi'')\\
&=&
\sum_{(\pi,\eta)\in {\cal DP}_{2}(n)}(-1)^{|\eta'|}m_{\mu}(\eta')m_{\nu}(\pi'')
\end{eqnarray*}
for every $n\geq 1$, where every $\eta$ has the decomposition
$\eta=\eta'\cup \eta''$ with every block of $\eta'$ obtained by connecting certain blocks of $\sigma'$
and $\eta''=\pi''$.

Let us finally demonstrate that $R_{n}=L_{n}$ for every $n\geq 1$.
In the expression for $R_{n}$ there is a summation over the set ${\cal F}(n)$
of triples $(m,\sigma , j)$, where $1\leq m \leq n$, $\sigma\in {\cal I}(m)$
and $j=(j_{1},j_{2}, \ldots , j_{m-1})$ is a tuple of non-negative integers
whose sum is equal to $n-m$. By Proposition 3.2, we get a bijection
$g: {\cal DP}_{2}(n) \rightarrow {\cal F}(n)$ which assigns
to every pair $(\pi,\eta)\in {\cal DP}_{2}(n)$ the triple $(m,\sigma,j)$.
Therefore, the summation over the set ${\cal F}(n)$ can be replaced by
a summation over the set ${\cal DP}_{2}(n)$ and since we can identify $\sigma$ with $\eta'$
as the proof of Proposition 3.2 shows, we have $m_{\mu}(\sigma)=m_{\mu}(\eta')$. Moreover,
$\nu(j_{1})\nu(j_{2})\ldots \nu(j_{m-1})=m_{\nu}(\pi'')$.
Therefore, $R_{n}=L_{n}$. This completes the proof.\hfill $\blacksquare$\\
\indent{\par}
{\sc Corollary 6.3.}
{\it In terms of $K$-transforms, the formula of Theorem 6.2 reads}
\begin{equation}
K_{\mu \,\vdash \,\nu}(z)=K_{\mu}(z-K_{\nu}(z))
\end{equation}
{\it where $K_{\mu}(z)$ and $K_{\nu}(z)$ are the $K$-transforms of $\mu$ and $\nu$, respectively.}\\
\indent{\par}
{\it Proof.}
This is an immediate consequence of Theorem 6.2 and the definition of the
$K$-transform. \hfill $\blacksquare$\\
\indent{\par}
In order to apply (6.2)-(6.3) to some examples, let
us recall basic facts on Jacobi continued fractions.
It is well-known [2] that
every probability measure $\mu$ with finite moments of all orders is
characterized by the sequences of Jacobi parameters $\alpha=(\alpha_{n})$ and
$\omega =(\omega_{n})$, $n\geq 0$, where $\alpha_{n}\in {\mathbb R}$ and
$\omega_{n}\geq 0$ (we will call them {\it Jacobi sequences}).
In that case we use the notation
$J(\mu)= (\alpha ,\omega)$. The Cauchy transform of $\mu$ can then
be expressed as a continued fraction of the form
\begin{equation}
G_{\mu}(z)=\cfrac{1}{z-\alpha_{0}-\cfrac{\omega_0}
{z-\alpha_{1}-\cfrac{\omega_{1}}{z-\alpha_{2}-\cfrac{\omega_{2}} {\ldots}}}}
\end{equation}
and it is understood that if $\omega_{m}=0$ for some $m$, then the fraction
terminates and, for convenience, we set $\omega_{n}=\alpha_{n}=0$ for all $n>m$.
In examples, we will mainly characterize measures by giving their Jacobi sequences
and refer the reader to [13] for details and explicit measures.

Let us also introduce the finite approximations of continued fractions (we shall use them
in Section 7 and 8).
In the case of the Cauchy transform $G_{\mu}(z)$  we define them as
quotients of polynomials of the form
\begin{equation}
[G_{\mu}(z)]_{m}=\frac{N_{m}(z)}{M_{m}(z)}, \;\; m\geq 1,
\end{equation}
where the numerators and the denominators satisfy the same recurrence
$$
Y_{m+1}(z)=(z-\alpha_{m})Y_{k}-\omega_{m-1}Y_{m-1}, \;\; m\geq 1,
$$
with different initial conditions: $N_{0}(z)=0, N_{1}(z)=1$
and $M_{0}(z)=1$, $M_{1}(z)=z-\alpha_{0}$ (see [2]).
In a similar way we define approximations of arbitrary continued fractions
and expressions involving them. In particular, we have
\begin{equation}
[F_{\mu}(z)]_{m-1}=z-[K_{\mu}(z)]_{m-1}=\frac{1}{[G_{\mu}(z)]_{m}}, \;\; m \geq 1,
\end{equation}
for the approximations of $F_{\mu}(z)$ and $K_{\mu}(z)$.

For any sequence $x=(x_{0},x_{1},x_{2}, \ldots)$ of real numbers,
let us also introduce the backward shift
$s(x)=(x_{1},x_{2},x_{3}, \ldots )$.
When we apply this shift to Jacobi sequences, we can find a relation 
between the orthogonal convolution of measures and the monotone convolution. \\
\indent{\par}
{\sc Corollary 6.4.}
{\it Let $\mu$ be a probability measure with finite moments of all orders
such that $J(\mu)=(\alpha,\omega)$. Then
$$
F_{\mu \, \vdash \,\nu}(z)=
z-\alpha_{0}-\frac{\omega_{0}}{F_{\mu_{s}\vartriangleright \nu}(z)}
$$
where $J(\mu_{s})=(s(\alpha), s(\omega))$, i.e. $\mu_{s}$
is a measure associated with shifted Jacobi sequences.}\\
\indent{\par}
{\it Proof.}
In (6.2), we write $F_{\mu}(w)$ as a continued fraction
and then substitute $w=F_{\nu}(z)$ to obtain
\begin{eqnarray*}
F_{\mu\,\vdash \,\nu}(z)
&=&
z-\alpha_{0}-\frac{\omega_{0}}{F_{\mu_{s}}(w)}\\
&=&
z-\alpha_{0}-\frac{\omega_{0}}{F_{\mu_{s}\vartriangleright \nu}(z)}
\end{eqnarray*}
which proves our assertion.  \hfill $\blacksquare$\\
\indent{\par}
{\sc Example 6.1}
Let $J(\mu)=(\alpha, \omega)$ and $\nu=\delta_{a}$,
where $a\in {\mathbb R}$. Then $F_{\nu}(z)=z-a$ and therefore, using Corollary 6.4, we get
$$
F_{\mu \,\vdash \, \nu}(z)=
z-\alpha_0-\frac{\omega_0}{F_{\mu_{s}}(z-a)}
$$
which shows that
$J(\mu\vdash\delta_{a})=((\alpha_{0}, \alpha_{1}+a, \alpha_{2}+a, \ldots ),\omega)$.
In particular, $\mu \vdash \delta_{0}=\mu$, i.e. $\delta_{0}$ is the
right identity w.r.t. the operation $\vdash$ (it is not hard to show that
the left identity does not exist). 
In turn, if $\mu=\delta_{a}$ and $J(\nu)=(\beta, \gamma)$,
where $a\in {\mathbb R}$, then $\alpha_{0}=a$ and $\omega_{0}=0$ and therefore, using Corollary 6.4,
we obtain $F_{\mu\,\vdash\,\nu}(z)=z-a$, which gives $\delta_{a}\vdash \nu=\delta_{a}$. \\
\indent{\par}
{\sc Example 6.2.}
Let $\mu = p\delta_{\lambda_{1}}+q\delta_{\lambda_{2}}$,
where $p+q=1$ and $J(\nu)=(\beta, \gamma)$. In that case
$$
G_{\mu}(z)=\frac{p}{z-\lambda_{1}}+\frac{q}{z-\lambda_{2}}
$$
and the reciprocal Cauchy transform is of the form
$$
F_{\mu}(z)=z-\lambda_{1}p-\lambda_{2}q-\frac{pq(\lambda_{1}-\lambda_{2})^{2}}{z-\lambda_{1}q-\lambda_{2}p}.
$$
Using (6.2), we obtain
$$
F_{\mu \,\vdash\, \nu}(z)=z-\lambda_{1}p-\lambda_{2}q-\frac{pq(\lambda_{1}-\lambda_{2})^{2}}
{F_{\nu}(z)-\lambda_{1}q-\lambda_{2}p}.
$$
Therefore,
$$
J(\mu \vdash \nu )=
\left((\lambda_{1}p+\lambda_{2}q, \beta_{0}+\lambda_{1}q+\lambda_{2}p, \beta_{1}, \beta_{2}, \ldots),
(pq(\lambda_{1}-\lambda_{2})^{2}, \gamma_{0}, \gamma_{1}, \ldots)\right).
$$
In particular, if $F_{\mu}(z)=z-\alpha_{0}-\omega_{0}/z$, then
$$
F_{\mu \vdash \nu}(z)=z-\alpha_{0}-\frac{\omega_{0}}{F_{\nu}(z)}
$$
and thus $J(\mu \,\vdash \,\nu)= ((\alpha_{0},\beta_{0}, \beta_{1}, \ldots),
(\omega_{0}, \gamma_{0}, \gamma_{1}, \ldots ))$.\\
\indent{\par}
{\sc Example 6.3.}
A closer look at equation (6.2) shows that it is natural to consider the orthogonal
convolution of measures which correspond to {\it mixed periodic J-fractions} 
[13] related to each other as follows:
\begin{eqnarray*}
J(\mu) &=& ((\alpha_{0}, \alpha_{1}, \alpha, \alpha, \ldots), (\omega_{0}, \omega_{1}, \omega, \omega, \ldots))\\
J(\nu) &=& ((\beta, \beta+\alpha, \beta+\alpha, \ldots ), (\gamma, \gamma+\omega, \gamma+\omega, \ldots))
\end{eqnarray*}
(here, $\alpha, \omega, \beta , \gamma$ denote numbers, not sequences).
In that case we have
\begin{eqnarray*}
F_{\mu}(z)&=&z-\alpha_{0}-\frac{\omega_{0}}{z-\alpha_{1}-\omega_{1} W_{(\alpha, \omega)}(z)}\\
F_{\nu}(z)&=&z-\beta-\gamma W_{(\alpha+\beta, \omega+\gamma)}(z)
\end{eqnarray*}
where $W_{(a,b)}(z)$ denotes the Cauchy transform of the Wigner measure $\sigma$
with mean $a$ and variance $b$.
In this case $J(\mu_{s})= ((\alpha_{1}, \alpha,\alpha, \ldots), (\omega_{1},\omega, \omega, \ldots))$ and
thus
\begin{eqnarray*}
F_{\mu_{s}\vartriangleright \nu}(z)
&=&
F_{\nu}(z)-\alpha_{1}-\omega_{1}W_{(\alpha,\omega)}(F_{\nu}(z))\\
&=&
z-\beta-\gamma W_{(\alpha+\beta, \omega+\gamma)}(z)-\alpha_{1}-\omega_{1}W_{(\alpha,\omega)}(F_{\nu}(z))\\
&=&
z-\beta-\alpha_{1}-(\gamma+\omega_{1}) W_{(\alpha+\beta, \omega+\gamma)}(z)
\end{eqnarray*}
since $W_{(\alpha+\beta, \omega+\gamma)}(z)=W_{(\alpha, \omega)}(F_{\nu}(z))$. Therefore, we obtain another
mixed periodic J-fraction
$J(\mu\vdash \nu) =((\alpha_{0},\alpha_{1}+\beta, \alpha+\beta, \alpha+\beta, \ldots ), (\omega_{0}, \omega_{1}+\gamma,
\omega+\gamma, \omega +\gamma, \ldots ))$.
For a discussion on the corresponding  measures, see [13].\\
\indent{\par}
The notion of the orthogonal convolution can be extended to the class of
all probability measures. Namely, by $\mu \vdash \nu$ we then understand the unique probability
measure defined by the reciprocal Cauchy transform of the form (6.2) - that the formula (6.2)
gives in fact a function from class ${\cal RC}$ is proven below. Note that the binary operation
$\vdash$ is neither commutative nor associative.\\
\indent{\par}
{\sc Theorem 6.5.} {\it If $\mu$ and $\nu$ are probability measures on the
real line, then the function of the form}
\begin{equation}
F(z)=F_{\mu}(F_{\nu}(z))-F_{\nu}(z)+z
\end{equation}
{\it defined on ${\mathbb C}^{+}$, is the reciprocal Cauchy transform of a probability
measure on the real line.}\\
\indent{\par}
{\it Proof.} In order to
demonstrate that $F(z)$ is the reciprocal of the Cauchy transform of a
probability measure, we will use the sufficiency condition of Maassen [18]
and show that
$$
\inf_{z\in {\mathbb C}^{+}}\frac{\Im(F(z))}{\Im z}=1
$$
where $\Im (u)$ denotes the imaginary part of $u\in {\mathbb C}$. Denoting
$w=F_{\nu}(z)$ and using the Nevanlinna representation theorem, we can write $F(z)$ in the
form
\begin{eqnarray*} F(z)&=&z-a-\int_{{\mathbb R}}\frac{1+xw}{w-x}d\tau(x)\\ &=& z-a -\int_{{\mathbb R}} \frac
{(1+xw)(\bar{w}-x)} {(w-x)(\bar{w}-x)} d\tau(x)
\end{eqnarray*}
where
$\tau$ is a positive finite measure, which gives $$ \Im F(z)=y+\Im
w\int_{{\mathbb R}}\frac{1+x^{2}}{|w-x|^{2}}d\tau(x)\geq y $$ for $z\in
{\mathbb C}^{+}$ since $\Im w=-\Im (G_{\nu}(z))/|G_{\nu}(z)|^{2}\geq 0$ (we
have $G_{\nu}: {\mathbb C}_{+}\rightarrow {\mathbb C}_{-}$). This implies
that
$$
\inf_{z\in {\mathbb C}^{+}}\frac{\Im F(z)}{\Im z}\geq 1.
$$
Moreover, we can write
$$
\frac{F(z)}{z}=1-\frac{F_{1}(z)}{z}+\frac{F_{2}(z)}{z}
$$
where
$$
F_{1}(z)=F_{\nu}(z) \;\; {\rm and} \;\; F_{2}(z)=F_{\mu}(F_{\nu}(z)) =
F_{\mu \vartriangleright\nu}(z)
$$
and thus both $F_{1}(z)$ and
$F_{2}(z)$ are reciprocals of Cauchy transforms of probability measures.
This gives
$$
\inf_{z\in {\mathbb C}^{+}}\frac{\Im F_{1}(z)}{\Im z}=1\;\;
{\rm and} \;\; \inf_{z\in {\mathbb C}^{+}}\frac{\Im F_{2}(z)}{\Im z}=1.
$$
Observe now that
$$
\Im (F_{1}(z)-F_{2}(z))=\Im
\left(K_{\mu}\left(\frac{1}{G_{\nu}(z)}\right)\right) \leq 0
$$
for $z\in
{\mathbb C}^{+}$ since $1/G_{\nu}:\; {\mathbb C}_{+}\rightarrow {\mathbb C}_{+}$
and $K_{\nu}:{\mathbb C}^{+}\rightarrow {\mathbb C}^{-}\cup
{\mathbb R}$. Therefore
$$
\frac{\Im(F_{2}(z)-F_{1}(z))}{\Im z}\geq 0,
\;\;\; z\in {\mathbb C}^{+}
$$
Now, since $F_{1}(z)$ and $F_{2}(z)$ are
holomorphic on ${\mathbb C}^{+}$, we have two real-valued functions,
$f_{1}(z):=\Im F_{1}(z)/\Im z$ and $f_{2}(z):=\Im F_{2}(z)/\Im z$, which are
continuous on ${\mathbb C}^{+}$ with $f_{2}\geq f_{1}$ on ${\mathbb C}^{+}$
and ${\rm inf}_{z\in {\mathbb C}^{+}}f_{1}(z)= {\rm inf}_{z\in {\mathbb
C}^{+}}f_{2}(z)=1$.
Thus, there exists a sequence $(z_{n})\subset {\mathbb C}^{+}$
such that $\lim_{n\rightarrow \infty}f_{2}(z_{n})=1$ and thus
$\lim_{n\rightarrow \infty}f_{1}(z_{n})=1$ which implies that
$\lim_{n\rightarrow \infty}f(z_{n})=0$, where $f=f_{2}-f_{1}$.
This proves hat
$$
\inf_{z\in {\mathbb C}^{+}}\frac{\Im F(z)}{\Im z}\leq 1.
$$
which completes the proof of the sufficiency condition. \hfill
$\blacksquare$\\
\indent{\par}
{\sc Corollary 6.6.}
{\it If $\mu$ and $\nu$ are probability measures on the real line, then the monotone convolution
of $\mu$ and $\nu$ can be decomposed as $\mu\vartriangleright\nu =
(\mu \vdash \nu)\uplus \nu$.}\\
\indent{\par}
{\it Proof.}
This decomposition is a direct consequence of (2.3),(2.5) and Theorem 6.5.
\hfill $\blacksquare$\\[10pt]
\newpage
\myownsection
\begin{center}
{\sc 7. Structures related to subordination functions}
\end{center}
In analogy to Sections 4 and 5, where we studied `orthogonal structures',
we now define and study structures (subalgebras, products and convolutions)
related to the subordination functions $F_{1}(z)$ and $F_{2}(z)$.
In particular, these functions uniquely determine probability measures which can 
be treated as convolutions of $\mu$ and $\nu$. 
These convolutions resemble the free additive convolution, except that one measure can be viewed as 
`subordinate' to the other. The same holds for the associated subalgebras and Hilbert spaces and this
motivates our terminology - `s-free convolution', `s-free subalgebras' and
`s-free product of Hilbert spaces'. In the case of compactly supported probability measures,
these structures can also be obtained as inductive limits of `alternating orthogonal structures',
but it seems to be of advantage to define them directly.

Let $({\cal H}_{i}, \xi_{i})$, $i=1,2$, be Hilbert spaces with distinguished unit vectors.
Then their Hilbert space free product $({\cal H}_{1}, \xi_{1})*({\cal H}_{2},\xi_{2})$
is $({\cal H},\xi)$ where
\begin{equation}
{\cal H}={\mathbb C}\;\xi\oplus\bigoplus_{n=1}^{\infty}
\bigoplus_{i_1 \ne i_2 \ne ... \ne i_n}
{\cal H}_{i_1}^0\otimes {\cal H}_{i_2}^0\otimes \ldots \otimes {\cal H}_{i_n}^0
\end{equation}
with ${\cal H}_i^0={\cal H}_i \ominus {\mathbb C} \xi_i$ and $\xi$ denoting a unit vector
(canonical scalar product is used).
For any $h\in{\cal H}_i$, denote by $h^0$ the orthogonal projection of $h$ onto ${\cal H}_i^0$.
Moreover, let
\begin{equation}
{\cal H}^{(n)}(j)=
\bigoplus_{\stackrel {i_{1}\neq i_{2}\neq \ldots \neq i_{n}}
{\scriptscriptstyle i_{1}\neq j}}
{\cal H}_{i_1}^0\otimes {\cal H}_{i_2}^0\otimes \ldots \otimes {\cal H}_{i_n}^0
\end{equation}
\begin{equation}
{\cal K}^{(n)}(j)=
\bigoplus_{\stackrel {i_{1}\neq i_{2}\neq \ldots \neq i_{n}}
{\scriptscriptstyle i_{n}\neq j}}
{\cal H}_{i_1}^0\otimes {\cal H}_{i_2}^0\otimes \ldots \otimes {\cal H}_{i_n}^0
\end{equation}
for any $j=1,2$, and $m\in {\mathbb N}$. For convenience,
also set ${\cal H}^{(0)}(j)={\cal K}^{(0)}(j)={\mathbb C}\xi$ with
the canonical projection $P_{0}:{\cal H}\rightarrow {\mathbb C}\xi$.
We will also use
\begin{equation}
{\cal H}(j)=\bigoplus_{n=1}^{\infty}{\cal H}^{(n)}(j), \;\;\;
{\cal K}(j)=\bigoplus_{n=1}^{\infty}{\cal K}^{(n)}(j)
\end{equation}
for $j=1,2$, and ${\cal H}^{(n)}={\cal H}^{(n)}(1)\oplus {\cal H}^{(n)}(2)$ for
$n\in {\mathbb N}$. Thus, we have
$$
{\cal H}=
{\mathbb C}\xi\oplus {\cal H}(1)\oplus {\cal H}(2)=
{\mathbb C}\xi\oplus {\cal K}(1)\oplus {\cal K}(2) 
$$
hence ${\cal H}$ can be decomposed as the union of ${\mathbb C}\xi$
and two `branches' (originating or ending with ${\cal H}_{1}^{0}$ or ${\cal H}_{2}^{0}$).

Using these notations, we can also decompose the free product of Hilbert spaces (7.1) as the
orthogonal direct sums
\begin{equation}
{\cal H}=\bigoplus_{n\geq 1}{\cal H}^{(n-1)}(j)\oplus {\cal H}^{(n)}(\bar{j})
\end{equation}
where we adopt the notation $\bar{1}=2$ and $\bar{2}=1$.
Moreover, let $P_{j}(n)$ denote the orthogonal projection onto
${\cal H}^{(n-1)}(j)\oplus {\cal H}^{(n)}(\bar{j})$.

The above notations will always be used in the following context.
Let $(\mathcal{A}_i,\varphi_i)$, $i=1,2$, be $\mathcal{C}^*$-noncommutative probability spaces.
We denote by $(\mathcal{H}_i,\pi_i,\xi_i)$ the GNS triple
of $(\mathcal{A}_i, \varphi_{i})$, i.e. $\mathcal{H}_i$
is a Hilbert space, $\xi_i$ is a cyclic (unit) vector in
$\mathcal{H}_i$ and $\pi_i:\mathcal{A}_i\rightarrow\mathcal{B}(\mathcal{H}_i)$ is a *-homomorphism, such that
$\varphi_i(a)=\langle \pi_i(a)\xi_i , \xi_i \rangle$ for all $a\in\mathcal{A}_i$, where
$\langle\cdot\,,\cdot\rangle$ denotes the scalar product in $\mathcal{H}_i$ (for simplicity, we use
the same notation for all scalar products).

Recall the definition of the free product representation.
On ${\cal H}$ we define a *-representation $\lambda_i:\mathcal{A}_i\rightarrow\mathcal{B}(\mathcal{H})$
of each algebra $\mathcal{A}_i$, $i=1,2$, as follows:
\begin{eqnarray*}
\lambda_{i}(a)(h_{1}\otimes h)&=&(\pi_{i}(a)h_{1})^{0} \otimes h + \langle \pi_{i}(a)h_{i},\xi_{i}\rangle h\\
\lambda_{i}(a)(\xi)&=&(\pi_{i}(a)\xi_{i})^{0}  + \langle \pi_{i}(a)\xi_{i},\xi_{i}\rangle \xi
\end{eqnarray*}
for any $h\in {\cal H}(i)$ and $h_{1}\in {\cal H}_{i}$, where identifcations
$\xi_{i}\otimes h\equiv h$ are made for any
$h_{1}\in {\cal H}_{i}^{0}$ and $h\in {\cal H}(i)$.
The {\it free product} of $(\lambda_{i})_{i\in I}$ is the representation
$\lambda=*_{i\in I}\lambda_{i}:*_{i\in I}{\cal A}_{i} \rightarrow {\cal B}({\cal H})$ given by
the linear extension of
$$
(*_{i\in I}\lambda_{i})(a_{1}a_{2}\ldots a_{n})=
\lambda_{i_{1}}(a_{1})\lambda_{i_{2}}(a_{2})\ldots \lambda_{i_{n}}(a_{n})
$$
for $a_{j}\in {\cal A}_{i_{j}}$, where $i_{1}, i_{2}, \ldots , i_{n}$ is a sequence of
alternating $1$'s and $2$'s.
Finally, on $\mathcal{B}(\mathcal{H})$ we define the so-called {\it vacuum state}
$\varphi(\cdot)=\langle \cdot\Omega,\Omega\rangle$ and the {\it free product}
of states $(\varphi_{i})_{i\in I}$ as the functional
$*_{i\in I}\varphi_{i}: *_{i\in I}{\cal A}_{i}\rightarrow {\mathbb C}$
given by the composition $ *_{i\in I}\varphi_{i}=\varphi \circ \lambda$.

Let $a_1\in {\cal B}({\cal H}_{1})$ and $a_2\in {\cal B}({\cal H}_{2})$ be fixed
random variables with distributions $\mu$ and $\nu$, respectively.
The corresponding free random variables $\lambda(a_1)$ and $\lambda(a_2)$
can be decomposed according to the Hilbert space decompositions (7.5) and
can be interpreted as consisting of sums of replicas of $a_{1}$ and $a_{2}$,
respectively (see Proposition 7.1). 
\\
\indent{\par}
{\sc Proposition 7.1.}
{\it According to the decomposition (7.5), the free random variable $\lambda(a_j)$ is the strongly
convergent series}
\begin{equation}
\lambda(a_j)=\sum_{n=1}^{\infty} a_j(n)
\end{equation}
{\it where $a_j(n) = P_{j}(n)a_jP_{j}(n)$ are replicas of $a_{j}$, where $j=1,2$.}\\
\indent{\par}
{\it Proof.}
Note that the subspace $\,{\cal H}^{(n-1)}(j)\oplus {\cal H}^{(n)}(\bar{j})\,$ is left invariant
by $\lambda(a_j)$ for every $n\in {\mathbb N}$. Using the direct sum decomposition (7.5), we get
the assertion.
\hfill $\blacksquare$\\
\indent{\par}
{\sc Remark 7.1.}
Decompositions of free random variables of type given by (7.6) were studied in the algebraic
framework of *-algebras [15], where it was shown that they can be viewed as `closed operators'
w.r.t. a suitable topology implemented by a sequence of projections.\\
\indent{\par}
From now on, when speaking of free random variables $a_{1},a_{2}$
as elements of ${\cal B}({\cal H})$, we will
understand that $a_j\equiv \lambda(a_j)\in {\cal B}({\cal H})$, where $j=1,2$. The Hilbert space
setting will be used below to study the `free convolution product' of $a_{1}$ and $a_{2}$ as well
as the `subordinate convolution' related to the subordination functions. We begin, however, with
a general formulation. \\
\indent{\par}
{\sc Definition 7.1.}
Let $({\cal A},\varphi, \psi)$ be a unital algebra with a pair
of linear normalized functionals.
Let ${\cal A}_{1}$ be a unital subalgebra of ${\cal A}$ and let
${\cal A}_{2}$ be a non-unital subalgebra with an `internal' unit $1_{2}$, i.e.
$1_{2}b=b=b1_{2}$ for every $b\in {\cal A}_{2}$.
We say that the pair $({\cal A}_{1},{\cal A}_{2})$ is {\it free with subordination}, or 
simply {\it s-free},
with respect to $(\varphi, \psi)$ if $\psi(1_{2})=1$ and it holds that\\
\indent{\par}
(i)
$\varphi(a_{1}a_{2}\ldots a_{n})=0$ whenever $a_{j}\in {\cal A}_{i_{j}}^{0}$ and
$i_{1}\ne i_{2}\neq \ldots \neq i_{n}$
\indent{\par}
(ii)
$\varphi(w_{1}1_{2}w_{2})=\varphi(w_{1}w_{2})-\varphi(w_{1})\varphi(w_{2})$
for any $w_{1},w_{2}\in {\rm alg}({\cal A}_{1}, {\cal A}_{2})$,\\[5pt]
where ${\cal A}_{1}^{0}={\cal A}_{1}\cap {\rm ker}\varphi$ and ${\cal A}_{2}^{0}={\cal A}_{2}\cap {\rm ker}\psi$.
We say that the pair $(a,b)$ of random variables from ${\cal A}$
is {\it s-free} with respect to $(\varphi, \psi)$
if the unital algebra generated by $a$ and the (non-unital) algebra 
generated by $b$ have this property.\\
\indent{\par}
The notion of `s-freeness' reminds freeness except that the `internal' unit $1_{2}$ 
is mapped by the GNS representation onto $P_{\xi}^{\perp}=1-P_{\xi}$, where $\xi$ is the distinguished 
unit vector of the Hilbert space, instead of the unit $1$ (see below).
Note also that condition (ii) resembles condition (ii) of Definition 4.1, but it is weaker
since it has $1_{2}$ `in the middle' and not an arbitrary $b\in{\cal A}_{2}$.
In particular, since ${\cal A}_{1}$ is unital, it also follows from (ii) that
$\varphi$ vanishes on ${\cal A}_{2}$ (cf. (ii) of Definition 4.1). 
Let us also point out that by conditions (i)-(ii) of Definition 7.1, $\varphi$ is uniquely determined on 
${\rm alg}({\cal A}_{1},{\cal A}_{2})$
by restrictions $\varphi|{\cal A}_{1}$ and $\psi|{\cal A}_{2}$ and it vanishes on ${\cal A}_{2}$
(as in the orthogonal case). 
Finally, note that `s-freeness' w.r.t. $(\varphi, \psi)$ differs from `conditional 
freeness' w.r.t $(\varphi, \psi)$, as in the orthogonal case (in particular, it is not symmetric w.r.t.
${\cal A}_{1}$ and ${\cal A}_{2}$).

The corresponding Hilbert space setting can be given as follows.\\
\indent{\par}
{\sc Definition 7.2.}
Let $({\cal H}_{1}, \xi_{1})$ and $({\cal H}_{2}, \xi_{2})$ be
Hilbert spaces with distinguished unit vectors $\xi_{1}$ and $\xi_{2}$, respectively.
The {\it s-free product} of $({\cal H}_{1}, \xi_{1})$ and $({\cal H}_{2}, \xi_{2})$
is the pair $({\cal K}, \xi)$, where ${\cal K}={\mathbb C}\xi \oplus {\cal K}(2)$.
We denote it by $({\cal K}, \xi)=({\cal H}_{1}, \xi_{1}) \oright ({\cal H}_{2}, \xi_{2})$
and by $\varphi$ - the canonical state on ${\cal B}({\cal H})$ associated with $\xi$.\\
\indent{\par}
Let us define *-representations
$\rho_{i}:{\cal B}({\cal H}_{i})\rightarrow {\cal B}({\cal K})$
by strongly convergent series
\begin{equation}
\rho_{1}(a_1)= \sum_{r=0}^{\infty}a_1(2r+1), \;\;\;
\rho_{2}(a_2)= \sum_{r=1}^{\infty}a_2(2r)
\end{equation}
where $a_1\in {\cal B}({\cal H}_{1})$, $a_2\in {\cal B}({\cal H}_{2})$.
Note that $\rho_{1}$ ($\rho_{2}$) is a faithful unital (non-unital)
*-homomorphism. Using these representations, we can describe the properties of the
s-free product of Hilbert spaces. Of course, if we consider
$({\cal H}_{2}, \xi_{2})\oright ({\cal H}_{1}, \xi_{1})$, we need to define
a different pair of *-representations (with the roles of ${\cal B}({\cal H}_{1})$ and
${\cal B}({\cal H}_{2})$ interchanged).\\
\indent{\par}
{\sc Theorem 7.2.}
{\it Let $({\cal K}, \xi)= ({\cal H}_{1},\xi_{1})\oright ({\cal H}_{2}, \xi_{2})$ be the
s-free product of Hilbert spaces and let $\varphi$, $\varphi_{1}$ and $\varphi_{2}$ be
the states associated with unit vectors $\xi$, $\xi_{1}$ and $\xi_{2}$, respectively.
Moreover, let $\psi$ be the state on ${\cal B}({\cal K})$
associated with any unit vector $\eta\in {\cal H}_{1}^{0}\subset {\cal H}$.
Finally, let ${\cal A}_{i}=\rho_{i}({\cal B}({\cal H}_{i}))$, where $i=1,2$.
Then
\indent{\par}
(i) \;\;the pair $({\cal A}_{1},{\cal A}_{2})$ is s-free w.r.t. $(\varphi, \psi)$,
\indent{\par}
(ii) \;$\varphi \circ \rho_{1}$ agrees with the expectation $\varphi_{1}$ on ${\cal B}({\cal H}_{1})$,
\indent{\par}
(iii) $\psi \circ \rho_{2}$ agrees with the expectation $\varphi_{2}$ on ${\cal B}({\cal H}_{2})$.}\\
\indent{\par}
{\it Proof.}
Let $a_{j}\in {\cal A}_{i_{j}}^{0}$, $1\leq j \leq n$ with $i_{1}\neq i_{2}\neq \ldots \neq i_{n}$.
In order to see that condition (i) of Definition 7.1 holds, first observe that
it holds whenever $a_{n}\in {\cal A}_{2}^{0}$,
since $\rho_{2}(b)\xi=0$ for any $b\in {\cal B}({\cal H}_{2})$.
Therefore, assume that $a_{n}\in {\cal A}_{1}^{0}$. In that case
$a_{n}\xi=h_{n}\in {\cal H}_{1}^{0}$,
$a_{n-1}h_{n}=h_{n-1}\otimes h_{n}\in {\cal H}_{2}^{0}\otimes {\cal H}_{1}^{0}$, $\;a_{n-2}h_{n-1}\otimes h_{n}=
h_{n-2}\otimes h_{n-1}\otimes h_{n}\in {\cal H}_{1}^{0}\otimes {\cal H}_{2}^{0}\otimes {\cal H}_{1}^{0}$, etc.
Continuing this process, we get condition (i) of Definition 7.1.
Finally, let $w_{1},w_{2}\in {\rm alg}({\cal A}_{1},{\cal A}_{2})$ and
observe that $1_{2}=\rho_{2}(1_{{\cal H}_{2}})=1-P_{\xi}$. This proves condition (ii) of Definition 7.1
and thus completes the proof of (i).
Verification of (ii) and (iii) is straightforward.
\hfill $\blacksquare$\\
\indent{\par}
{\sc Corollary 7.3.}
{\it Let $\mu, \nu$ be compactly supported probability measures on ${\mathbb R}$.
Then there exist a Hilbert space ${\cal K}$, unit vectors $\xi, \eta\in {\cal K}$ and
self-adjoint bounded random variables $X_1,X_2\in {\cal B}({\cal K})$ such that the pair $(X_1,X_2)$ is s-free
w.r.t. $(\varphi, \psi)$, where $\varphi$ and $\psi$ are vector states associated with
$\xi, \eta\in {\cal K}$. Moroever, the $\varphi$-distribution of $X_1$ and
the $\psi$-distribution of $X_2$ coincide with $\mu$ and $\nu$, respectively. Finally,
the $\varphi$-distribution of $X_1+X_2$, denoted $\mu\, \boxright \,\nu$, is compactly supported.}\\
\indent{\par}
{\it Proof.}
The proof is similar to that of Corollary 4.2 - replace ${\cal H}$ by 
${\cal K}={\mathbb C}\xi \oplus {\cal K}(2)$ and
*-homomorphisms $\tau_{1}$ and $\tau_{2}$ by $\rho_{1}$ and $\rho_{2}$, respectively.
\hfill $\blacksquare$\\
\indent{\par}
The s-free product of $C^{*}$-probability spaces can be defined 
along the lines of Section 4 (together with a $C^{*}$-version of Theorem 7.2). 
Computations of convolutions $\mu\,\boxright \,\nu$ are postponed till Section 8,
where the transforms are studied, which provide the natural tools.
However, we give here an example which shows a natural connection between the
s-free product and branches of the free product of graphs (convolutions
$\mu\,\boxright \,\nu$ will then give their spectral distributions). \\
\indent{\par}
{\sc Example 7.1.}
Consider two rooted graphs as in Example 4.1. Let $({\cal B}_{1},e)$
be the (uniformly locally finite) rooted graph called `branch subordinate to ${\cal G}_{1}$'
(or, simply, `branch') obtained from the free product of rooted graphs
$({\cal G}_{1},e_{1})*({\cal G}_{2},e_{2})$ by restricting the set of vertices $V_{1}*V_{2}$
to the set $V$ consisting of the empty word $e$ and words ending
with a vertex (letter) from $V_{1}^{0}$. Then the adjacency matrix of $({\cal B}_{1},e)$
can be decomposed as $A({\cal B}_{1})=A^{(1)}+A^{(2)}$, where the summands
$$
A^{(1)}=\sum_{n\; {\rm odd}}A_{1}(n), \;\;\;
A^{(2)}=\sum_{n\; {\rm even}}A_{2}(n),
$$
(with $n$ positive) are bounded operators on $l_{2}(V)$ which are s-free w.r.t.
$(\varphi, \psi)$, where $\varphi(.)=\langle .\delta (e),\delta (e)\rangle$ and 
$\psi(.)=\langle . \delta (v),\delta (v)\rangle$ with $v\in V_{1}^{0}$. 
Moroever, the pair $({\cal A}, \varphi)$, where ${\cal A}$ is the $C^{*}$-algebra generated by 
$A^{(1)}$, $A^{(2)}$ and the identity $I$ on $l_{2}(V)$, is the s-free product of the $C^{*}$-probability 
spaces $({\cal A}_{1},\varphi_{1})$ and $({\cal A}_{2}, \varphi_{2})$, where 
${\cal A}_{i}$ is generated by $A_{i}$ and the unit $I_{i}$ on $l_{2}(V_{i})$,
and $\varphi_{i}$ is the state defined by the vector $\delta (e_{i})$, $i=1,2$.
Thus the branch ${\cal B}_{1}$ can be viewed as the s-free product of two subgraphs
obtained from the coverings of ${\cal B}_{1}$ built from replicas of ${\cal G}_{1}$ and ${\cal G}_{2}$,
respectively. The `branch subordinate to ${\cal G}_{2}$' can be decomposed in a similar way. 
A detailed study of `branches' of the free product of rooted graphs will be given in [1].\\[10pt]
\myownsection
\begin{center}
{\sc 8. Subordinate branches}
\end{center}
\indent{\par}
In this Section we study `complete' decompositions of the sum $X_{1}+X_{2}$ of
free bounded random variables with distributions $\mu$ and $\nu$, respectively.
This leads to the concept of `free branches subordinate to $X_{1}$ or $X_{2}$'
(terminology is motivated by Example 7.1) which give a Hilbert space realization of the decomposition
of the s-free convolution of $\mu$ and $\nu$ into a sequence of orthogonal 
convolutions of alternating $\mu$ and $\nu$.\\
\indent{\par}
{\sc Definition 8.1.}
Let $(X_{1}, X_{2})$ be a pair of self-adjoint random variables from ${\cal B}({\cal H})$ which are
free w.r.t. $\varphi$. The self-adjoint random variable from ${\cal B}({\cal H})$ given by
the strongly convergent series
\begin{equation}
B_{j}(k)=\sum_{r=0}^{\infty}X_j(2r+k)+\sum_{r=0}^{\infty}X_{\bar{j}}(2r+k+1)
\end{equation}
where $j=1,2$ and $k\in {\mathbb N}$, is
called the {\it $k$-th free branch subordinate to $X_{j}$}.\\
\indent{\par}
In particular, for $k=1$, the 1-st free branch subordinate to $X_{1}$ takes the form
$$
B_{1}(1)=\rho_{1}(X_{1})+\rho_{2}(X_{2})
$$
which, for simplicity, we shall denote $B_{1}$.
Free branches of higher orders can then be obtained from the recursions
$$
B_{j}(k) = X_j(k) + B_{\bar{j}}(k+1)
$$
which shows that we get two disjoint sequences of `free branches' with alternating
subordination, begining either with $B_{1}\equiv B_{1}(1)$ or $B_{2}\equiv B_{2}(1)$.

Viewing the `free convolution product' of $X_{1}$ and $X_{2}$ as a tree-like structure
with two types of `leaves' (which it really is when $X_{1}$ and $X_{2}$ are taken to be 
adjacency matrices of graphs, cf. Example 7.1),
one can interpret free branches as follows.
The branch $B_{1}$ is the `half' of the `tree' which
originates with `leaf' $X_1(1)$ (the first replica of $X_1$),
followed by `leaves' produced by $X_{2}(2)$ (the second replica of $X_2$),
then `leaves' produced by $X_1(3)$ (the third replica of $X_3$), etc.
In the case of branches of higher order,
$B_{1}(k)$ is the part of the `tree' which originates from `leaves' $X_1(k)$ at height $k$,
followed by `leaves' $X_{2}(k+1)$, then $X_1(k+2)$, etc (it is similar for branches
which are `subordinate' to $X_{2}$).\\
\indent{\par}
{\sc Proposition 8.1.}
{\it Let $(X_{1},X_{2})$ be a pair of free
self-adjoint bounded random variables from ${\cal B}({\cal H})$
with distributions $\mu$ and $\nu$, respectively. For $k\in {\mathbb N}\cup \{0\}$,
let $\varphi_{k}$ and $\psi_{k}$ denote states on ${\cal B}({\cal H})$
associated with arbitrary unit vectors $\zeta_{k}\in {\cal H}^{(k)}(2)$ and $\eta_{k}\in {\cal H}^{(k)}(1)$,
respectively. Then the $\psi_{k-1}$-distribution of $B_{1}(k)$ and
the $\varphi_{k-1}$-distribution of $B_{2}(k)$ are given by
$\mu \, \boxright \,\nu$ and $\nu \, \boxright \, \mu$, respectively, for every $k\in {\mathbb N}$.} \\
\indent{\par}
{\it Proof.}
For $k=1$ the assertion is a direct consequence of Definition 8.1 and Corollary 7.3.
For $k>1$, without loss of generality, consider the $\psi_{k-1}$-distribution of
$B_{1}(k)$. We need to show that
$\psi_{k-1}(B_{1}^{n}(k))=\varphi (B_{1}^{n})$
for every $k\geq 2$ and $n\in {\mathbb N}$.
Therefore we can write
$$
B_{1}^{n}(k)\eta_{k-1}=(B_{1}^{n}\xi_{1}) \otimes \eta_{k-1}
$$
if we identify $\xi_{1}\otimes \eta_{k-1}$ with $\eta_{k-1}$. In fact, in the definition of $B_{1}(k)$,
we only have $X_1(j)$'s with $j\geq k$ and $X_2(j)$'s with $j\geq k+1$,
thus $X_1(k-1)$ is missing, which is the only summand among all
$X_1(r)$'s and $X_2(r)$'s in (7.6) which can map $\eta_{k-1}$ onto a vector which is not
in ${\mathbb C}\eta_{k-1}$. Therefore
$$
B_{1}^{n}(k)\eta_{k-1}=\varphi(B_{1}^{n})\eta_{k-1}\;\; {\rm mod}\;({\cal H}\ominus {\cal H}^{(k-1)}(1))
$$
which gives $\psi_{k-1}(B_{1}^{n}(k))=\varphi(B_{1}^{n})$ and that completes the proof.
\hfill $\blacksquare$\\
\indent{\par}
{\sc Lemma 8.2.}
{\it For every $k\in {\mathbb N}$, the pair of random variables $(X_1(k),B_{2}(k+1))$
is orthogonal with respect to $(\psi_{k-1}, \varphi_{k})$
and the pair $(X_2(k),B_{1}(k+1))$ is orthogonal with respect to
$(\varphi_{k-1}, \psi_{k})$.}\\
\indent{\par}
{\it Proof.}
Without loss of generality, consider the pair $(X_2(k),B_{1}(k+1))$.
For simplicity, denote $b=B_{1}(k+1)$, $y=Y(k)$
and identify each element $w$ of the algebra generated by $b$ and $y$ with $\lambda(w)$.
We need to show two orthogonality conditions of Definition 4.1:
$$
\varphi_{k-1}(w_1 b)=\varphi_{k-1}(bw_1)=0
$$
$$
\varphi_{k-1}(w_{1}yb^{n}yw_{2})=
\psi_{k}(b^{n})(\varphi_{k-1}(w_{1}y^2w_{2})-
\varphi_{k-1}(w_{1}y)\varphi_{k-1}(yw_{2}))
$$
for any $w_{1},w_{2}\in {\rm alg}(b,y)$. We have
$$
\varphi_{k-1}(w_{1}b)=\langle w_{1}b\zeta_{k-1}, \zeta_{k-1} \rangle
$$
where $\zeta_{k-1}\in {\cal H}^{(k-1)}(2)$ is a unit vector
(it is convenient to think of a simple tensor of lenght $k-1$ which begins with a vector
from ${\cal H}_{1}^{0}$).
By Definition 8.1, we have $b\zeta_{k-1}=0$.
This proves that $\varphi_{k-1}(w_{1}b)=0$ (the equation $\varphi_{k-1}(bw_{2})=0$ can
be obtained from this by taking the adjoints). This gives the first orthogonality condition.

To prove the second orthogonality condition, consider now $yw_{2}\zeta_{k-1}$.
Using the definition of $b$ again, we obtain
$$
yw_{2}\zeta_{k-1}=\alpha \zeta_{k-1}+ h_{2}\otimes \zeta_{k-1}
$$
for some $\alpha \in{\mathbb C}$ and  $h_{2}\in {\cal H}_{2}^{0}$.
Now, we know that $b^{n}\zeta_{k-1}=0$ for every $n\in {\mathbb N}$ since $b\zeta_{k-1}=0$,
and, moreover,
$$
b^{n}h_{2}\otimes \zeta_{k-1}=
\psi_{k}(b^{n})h_{2}\otimes \zeta_{k-1}\;\;\; {\rm mod}\; ({\cal H}^{(k+1)} \oplus \ldots \oplus {\cal H}^{(k+n)}).
$$
This gives
\begin{eqnarray*}
\varphi_{k-1}(w_{1}yb^{n}yw_{2})
&=&
\langle w_{1}yb^{n}yw_{2}\zeta_{k-1},\zeta_{k-1} \rangle\\
&=& \langle w_{1}yb^{n}(\alpha \zeta_{k-1} +h_{2}\otimes \zeta_{k-1}, \zeta_{k-1} \rangle\\
&=& \psi_{k}(b^{n})\langle w_{1}yh_{2}\otimes \zeta_{k-1},\zeta_{k-1} \rangle\\
&=&\psi_{k}(b^{n})\langle w_{1}y(yw_{2}\zeta_{k-1}- \alpha \zeta_{k-1}),\zeta_{k-1} \rangle\\
&=&
\psi_{k}(b^{n})
(\langle w_{1}y^2w_{2}\zeta_{k-1},\zeta_{k-1} \rangle -
\alpha
\langle w_{1}y\zeta_{k-1},\zeta_{k-1} \rangle)\\
&=&
\psi_{k}(b^{n})
(\varphi_{k-1}(w_{1}y^2w_{2})-\varphi_{k-1}(yw_{2})\varphi_{k-1}(w_{1}y))
\end{eqnarray*}
since $\alpha = \langle yw_{2}\zeta_{k-1},\zeta_{k-1}\rangle =\varphi_{k-1}(yw_{2})$. Thus, we proved 
the second orthogonality condition for $(X_2(k), B_{1}(k+1))$, which
completes the proof of the lemma.\hfill $\blacksquare$\\
\indent{\par}
{\sc Lemma 8.3.}
{\it Let $\mu$ and $\nu$ be the $\varphi$-distributions of $X_{1}$ and $X_{2}$, respectively, and
let $(\mu \vdash _{n}\nu)_{n\geq 1}$ be a sequence of distributions defined recursively by
\begin{equation}
\mu \vdash_{1}\nu =\mu \vdash \nu, \;\;\;
\mu \vdash_{n} \nu = \mu \vdash (\nu \vdash_{n-1} \mu),\;\;n\geq 2.
\end{equation}
Then the moments of $B_{1}$ and $B_{2}$ of orders $\leq 2m$ in the
state $\varphi$ agree with the corresponding moments of
$\mu \vdash_{m-1}\nu$ and $\nu \vdash_{m-1}\mu$, respectively.}\\
\indent{\par}
{\it Proof.}
Let us consider only the $\varphi$-distribution of $b=B_{1}$ since the proof
for $B_{2}$ is identical. We have
$$
b= Z_{1}+Z_{2}+ \ldots + Z_{m}
$$
where $Z_{1}=X_1(1), Z_{2}=X_2(2), Z_{3}=X_1(3), \ldots, Z_{m-1}=X_1(m-1), Z_{m}=B_{2}(m)$.
In view of Lemma 8.2, the pairs
$$
(Z_{1}, Z_{2}+\ldots + Z_{m}), \;(Z_{2}, Z_{3}+\ldots + Z_{m}), \ldots , (Z_{m-1},Z_{m})
$$
are orthogonal w.r.t.
$(\psi_{0}, \varphi_{1}), \;(\varphi_{1},\psi_{2}), \ldots , (\psi_{m-2}, \varphi_{m-1})$, respectively.
Thus,
$$
\varphi(b^{r})= \mu_{1}\vdash (\mu_{2} \vdash (\ldots (\mu_{m-1}\vdash \mu_{m})\ldots ))(r)
$$
for every natural $r$, where $\mu_{1},\mu_{2}, \ldots , \mu_{m}$ are the
$\psi_{0}$--, $\varphi_{1}$--, $\psi_{2}$--, $\ldots$, $\varphi_{m-1}$-- distributions of
$Z_{1}, Z_{2}, \ldots , Z_{m}$, respectively.
We know that $\mu_{1}=\mu$, $\mu_{2}=\nu$, $\ldots$, $\mu_{m-1}=\mu$
and $\mu_{m}$ is the distribution of $B_{2}(m)$. Thus, in view of Corollary 5.3, $\varphi(b^{r})$
depends on moments of $\mu_{i}$ of orders $1\leq j_{i} \leq r-2i+2$ for $1\leq i \leq [1/2(r+1)]$.
In particular, it depends on the moments of $\mu_{m}$ if and only if $r\geq 2m-1$.
However, if $r\leq 2m$, then only the moments of $\mu_{m}$ of orders $\leq 2$ may come into play.
The latter, however, agree with the moments of $\nu$. This proves our assertion.
\hfill $\blacksquare$\\
\indent{\par}
{\sc Theorem 8.4.}
{\it For compactly supported $\mu$ and $\nu$, it holds that}
\begin{equation}
K_{\mu\,\boxright\, \nu}(z)=
\lim_{m\rightarrow \infty}
K_{\mu\,\vdash_{m} \nu}(z)
\end{equation}
{\it where the convergence is uniform on compact subsets of ${\mathbb C}^{+}$.}\\
\indent{\par}
{\it Proof.}
In view of Lemma 8.3, we have convergence of moments and thus (since all measures involved have compact support) 
weak convergence
$$
w-\lim_{m\rightarrow \infty}(\mu \vdash _{m}\nu)= \mu \boxright \nu
$$
which implies that $K_{\mu \,\vdash_{m}\,\nu}(z)$ converges uniformly to $K_{\mu \,\boxplus \,\nu}(z)$
on compact subsets of ${\mathbb C}^{+}$. Using Corollary 6.3, we obtain
$$
K_{\mu\,\vdash_{m} \nu}(z)=
K_{\mu_1}(z-K_{\mu_2}(z-K_{\mu_3}(\ldots (K_{\mu_m}(z)))))
$$
where $\mu_{1}=\mu$, $\mu_{2}=\nu$, $\mu_{3}=\mu$, $\ldots$, $\mu_{m}=\mu$ (if $m$ even)
or $\mu_{m}=\nu$ (if $m$ odd). This completes the proof. \hfill $\blacksquare$\\
\indent{\par}
{\sc Remark 8.1.}
An informal way of writing (8.3) is to
use the `continued composition' form
\begin{equation}
K_{\mu\,\boxright\, \nu}(z)=
K_{\mu}(z-K_{\nu}(z-K_{\mu}(z-K_{\nu}(\ldots ))))
\end{equation}
which is particularly appealing when the $K$-transforms involved are simple
and can be easily associated with Jacobi continued fractions.\\
\indent{\par}
{\sc Remark 8.2.}
It follows from Lemma 8.3 that the moments of $\mu \,\boxright \,\nu$
depend only on the moments of $\mu$ and $\nu$.
Moreover, from Corollary 7.3 we know that $\boxright$ is a binary operation on ${\cal M}_{c}$
(neither commutative nor associative).
Finally, using Theorem 8.4 (see also Example 8.1 below),
we can see that $\delta_{0}$ is the right identity w.r.t. $\boxright$
(the left identity does not exist).\\
\indent{\par}
{\sc Remark 8.3.}
More generally, if $\mu$ and $\nu$ are arbitrary probability measures,
one can define $\mu \boxright \nu$ and $\nu \boxright \mu$
as the unique probability measures associated with the reciprocal Cauchy transforms
$F_{1}(z)$ and $F_{2}(z)$, respectively, and thus $\boxright$ extends to a binary operation on
${\cal M}$. Moreover, $\delta_{0}$ is the right identity w.r.t. the operation
$\boxright$ on all of ${\cal M}$ since $F_{1}(z)=z$, $F_{2}(z)=F_{\mu}(z)$ are the unique
functions from class ${\cal RC}$ which satisfy (1.8).\\
\indent{\par}
{\sc Example 8.1.}
If $\mu=\delta_{a}$ and $J(\nu)= (\beta, \gamma)$, then we can use properties
of the orthogonal convolution (see Example 6.1) to obtain
$\delta_{a}\,\boxright \,\nu=\delta_{a}\vdash (\nu \,\boxright \,\delta_{a})=\delta_{a}$ since
$\delta_{a}\vdash \sigma =\delta_{a}$ for any $\sigma$.
Similarly, if $\nu=\delta_{a}$ and $J(\mu)=(\alpha, \omega)$, we get
$\mu \,\boxright \,\delta_{a}=\mu \vdash (\delta_{a}\vdash (\mu \,\boxright \,\delta_{a}))=
\mu \vdash \delta_{a}$ and we already know that 
$J(\mu\vdash \delta_{a})=
((\alpha_{0},\alpha_{1}+a, \alpha_{2}+a, \ldots), \omega)$,
which gives $J(\mu \,\boxright\,\nu)$.\\
\indent{\par}
{\sc Example 8.2.}
Let $K_{\mu}(z)=\alpha+\omega/z$ and $K_{\nu}(z)=\beta +\gamma/z$, where $\omega\neq 0 \neq \gamma$.
Using Theorem 8.4 and the definition of the $K$-transform, we obtain
$$
K_{\mu \,\boxright \,\nu}(z)=
\alpha+\cfrac{\omega}{z-\beta-\cfrac{\gamma}{z-\alpha-\cfrac{\omega}{z-\beta-\cfrac{\gamma}{\ldots}}}}
$$
and therefore, $G_{\mu\, \boxright \, \nu}(z)$ is a 2-periodic J-fraction.
Algebraic calculations give
$$
G_{\mu\,\boxright\,\nu}(z)=\frac{P(z)-2\gamma-\sqrt{P^{2}(z)-4A^{2}}}{2\Delta(z)}
$$
where $P(z)=-\gamma-\omega +(z-\alpha)(z-\beta)$, $\Delta(z)=\gamma(z-\alpha)$ and $A^{2}=\gamma\omega$.
The absolutely continuous part of this measure is of the form
$$
f(x)=\frac{\sqrt{4A^{2}-P^{2}(x)}}{2\pi \gamma (x-\alpha)}
$$
on the `stable band' with end-points 
$\alpha+\beta\pm \sqrt{(\alpha-\beta)^{2}+4\gamma +4\omega +8A}$,
with a possible atom at $x=\alpha$ if $\gamma \neq \omega$. For details, see [13].\\
\indent{\par}
{\sc Example 8.3.}
Consider $\mu$ to be the Wigner measure with mean $\alpha$ and variance $\omega$.
Then it holds that
$J(\mu \boxright \mu)= ((\alpha, 2\alpha, 2\alpha, \ldots), (\omega, 2\omega, 2\omega, \ldots))$.
To show this, it is enough to derive a suitable formula for
$J_{m}(z):=[K_{\mu\,\boxright \,\mu}(z)]_{m}$.
Using the continued fraction
$$
K_{\mu}(z)=\alpha+\cfrac{\omega}{z-\alpha-\cfrac{\omega}{z-\alpha-\cfrac{\omega}{\ldots }}}
$$
and its finite approximations (6.6), we obtain by the induction argument
\begin{eqnarray*}
J_{m}(z)
&=&
[K_{\mu\,\vdash_{m-1}\,\mu}(z)]_{m}\\
&=&
[K_{\mu}(z-J_{m-1}(z))]_{m}\\
&=&
\alpha+\cfrac{\omega}{z-\alpha-J_{m-1}(z)-\cfrac{\omega}{z-\alpha-J_{m-2}(z)-\cfrac{\omega}
{\ldots}}}\\
&=&
\alpha+ \frac{\omega}{z-2J_{m-1}(z)}
\end{eqnarray*}
and thus $K_{\mu\,\boxright \,\mu}(z)=\omega W_{2\alpha,2\omega}(z)$,
which proves our claim that we get a mixed periodic J-fraction.
The analytic form of the Cauchy transform is
$$
G_{\mu\,\boxright \, \mu}(z)=
\frac{\Lambda_{1}(z)-\sqrt{P^{2}(z)-A^{2}}}
{2\Gamma_{1}(z)}
$$
where $P(z)=(z-2\alpha)^{2}$, $A^{2}=2\omega^{2}$, $\Lambda_{1}(z)=3z-2\alpha$
and $\Gamma_{1}(z)=z^{2}-z\alpha+\omega^{2}$, using the notation of [Theorem 3.4, Ref.13] 
(with $N=1$, $M=1$). Thus $\mu\,\boxright \,\mu$ is the measure with the absolutely continuous 
part given by the density
$$
f(x)=\frac
{\sqrt{A^{2}-P^{2}(x)}}
{2\pi \Gamma_{1}(x)}
$$
on the `stable band' with endpoints at $\alpha \pm 2\sqrt{2}\omega$,
with possible atoms at $1/2(\alpha \pm \sqrt{\alpha^{2}-4\omega^{2}})$ if $\alpha^{2}-4\omega^{2}\geq 0$.
\\
\indent{\par}
It follows from Theorem 8.4 that every compactly supported probability measure
whose Cauchy transform can be represented in the form of a 2-periodic continued fraction 
is the weak limit of a sequence $\mu\,\vdash_{m} \,\nu$ for suitably chosen $\mu$ and $\nu$.
This corollary can be easily generalized to make every compactly supported probability 
measure the weak limit of a sequence of $m$-fold orthogonal convolutions. 
It can be shown that this setting corresponds to the model
of freeness with infinitely many states [10], but we will not discuss this connection here.\\
\indent{\par}
{\sc Corollary 8.5.}
{\it Let $\mu$ be a compactly supported probability measure on the real line
such that $J(\mu)=(\alpha, \omega)$ and let $(\mu_{n})_{n\in {\mathbb N}}$ be 
discrete probability measures with $K$-transforms
$$
K_{\mu_1}(z)=\alpha_{0}+\frac{\omega_{0}}{z-\alpha_{1}}, \;\;\;
K_{\mu_{n}}(z)=\frac{\omega_{n-1}}{z-\alpha_{n}},\;\; n\geq 2.
$$
Then $\mu=w-\lim_{m\rightarrow \infty} \mu_{1}\vdash(\mu_{2}\vdash (\ldots \vdash \mu_{m}))$.}\\
\indent{\par}
{\it Proof.}
One can observe that the sequence of truncations of the $K$-transform  of $\mu$ has the form
$$
[K_{\mu}(z)]_{m}=K_{\mu_{1}}(z-K_{\mu_{2}}(z- \ldots (z-K_{\mu_m}(z))))=
K_{\mu_{1}\vdash(\mu_{2}\vdash (\ldots \vdash \mu_{m}))}(z)
$$
and therefore it converges uniformly to $K_{\mu}(z)$ on the compact subsets of ${\mathbb C}^{+}$,
from which the assertion follows.\hfill $\blacksquare$\\
\indent{\par}
{\sc Remark 8.4.}
We end this section with a comment on graphs and their products (cf. Examples 4.1,4.2,7.1)
which actually motivated some of our terminology.
In the context of rooted graphs, Theorem 8.4 corresponds to the 
`complete' decomposition of branches (briefly discussed in Example 7.1) 
as orthogonal products of (replicas of) alternating ${\cal G}_{1}$ 
and ${\cal G}_{2}$. Namely, taking two rooted graphs as in Example 4.1, we obtain
$$
{\cal B}_{1}=
{\cal G}_{1}\vdash ( {\cal G}_{2}\vdash({\cal G}_{1}\vdash ({\cal G}_{2}
\vdash (\ldots ))))
$$
in the case of ${\cal B}_{1}$.
The spectral distributions of this branch associated with the vector $\delta(e)$
is expressed in terms of $\mu$ and $\nu$ (with the notation of Example 4.1)
by their s-free convolution
$\mu \boxright \nu = w-\lim_{m\rightarrow \infty}(\mu \vdash _{m}\nu)$.
Similar formuals hold for ${\cal B}_{2}$.
For more details, see [1].\\[10pt]

\myownsection
\begin{center}
{\sc 9. Decomposition of the additive free convolution}
\end{center}
Using decompositions of free branches of Section 8, we can show now that the distribution of the sum $X_1+X_2$
of self-adjoint free random variables $X_1,X_2\in {\cal B}({\cal H})$ with
$\varphi$-distributions $\mu$ and $\nu$, respectively, can be `completely' decomposed.\\
\indent{\par}
{\sc Lemma 9.1.}
{\it In the decomposition $X_1+X_2=B_{1}+B_{2}$,
the branches $B_{1}$ and $B_{2}$ are boolean independent with respect to $\varphi$.
Therefore, $\mu \boxplus \nu = (\mu \boxright \nu) \uplus (\nu \boxright \mu)$ is the
corresponding decomposition of the free additive convolution.}\\
\indent{\par}
{\it Proof.}
Let $w$ be any element of the *-algebra generated by
$B_{1}$ and $B_{2}$ and, for simplicity, denote $x=X_{1}(1)$.
In order to show boolean independence of $B_{1}$ and $B_{2}$, we compute
the moment with $B_{1}^{k}$ at the `end' of the moment (the other case is similar):
\begin{eqnarray*}
\varphi(w B_2B_1^{k}) &=&
\langle w B_2B_{1}^{k}\xi , \xi\rangle \\
&=&
\langle w B_2x^{k}\xi , \xi\rangle\\
&=&
\langle x^{k}\xi , \xi\rangle  \langle w B_2\xi , \xi \rangle
+\langle w B_2(x^{k}\xi_{1})^{0}, \xi\rangle\\
&=&
\varphi(x^{k})
\varphi (w B_2)\\
&=&
\varphi(B_1^{k})
\varphi (w B_2)
\end{eqnarray*}
since $B_{2}h_{1}=0$ for every $h_{1}\in {\cal H}_{1}^{0}$. This completes the proof.
\hfill $\blacksquare$\\
\indent{\par}
{\sc Lemma 9.2.}
{\it In the decomposition $X_1+X_2=X_{j}(1)+Z_j$, where $j=1,2$,
the pair of random variables $(X_j(1),Z_j)$ is monotone independent w.r.t. $\varphi$. Therefore,
\begin{equation}
\mu \boxplus \nu = \mu \vartriangleright (\nu \boxright \mu), \;\;\;
\nu \boxplus \mu = \nu \vartriangleright (\mu \boxright \nu)
\end{equation}
are the corresponding decompositions of the free additive convolution for $j=1,2$, respectively.}\\
\indent{\par}
{\it Proof.}
It is enough to consider the decomposition for $j=1$ (the proof for $j=2$ is similar).
For simplicity, denote $z=Z_{1}$, $x=X_{1}(1)$, $b_{1}=B_{2}=B_{2}(1)$, $b_{2}=B_{2}(2)$
and $p=P_{1}(1)$. In order to show that (2.2) holds for the pair $(x,z)$, consider two cases.
\indent{\par}
{\it Case 1.} Suppose that $z^{n}$, $n\in {\mathbb N}$, is at the `end' of the moment. We need to show that
\begin{equation}
\varphi (wxz^{n})=\varphi (wx)\varphi (z^{n})
\end{equation}
for any $w\in {\rm alg}(x,z)$.
Write
$$
\varphi (wxz^{n})=\langle wxz^{n}\xi , \xi \rangle
$$
and examine $z^{n}\xi$.
Since $z=b_{1}+b_{2}$ does not contain $X_{1}(1)$, it holds that
$$
z^{n}\xi = b_{1}^{n}\xi\in {\mathbb C}\xi \oplus {\cal K}^{(1)}(1)\oplus \ldots \oplus {\cal K}^{(n)}(1)
$$
which implies that
\begin{equation}
z^{n}\xi=\varphi(z^{n})\xi \;\;  {\rm mod}\,({\cal K}^{(1)}(1)\oplus \ldots \oplus {\cal K}^{(n)}(1))
\end{equation}
for any $n\in {\mathbb N}$. In particular, this gives
$\varphi(z^{n})=\varphi(b_{1}^{n})$, i.e. that $z$ and $b_{1}$ are $\varphi$-identically distributed.
Now, since $x=pXp$ and $p:{\cal H}\rightarrow {\mathbb C}\xi \oplus {\cal H}_{1}^{0}$, we have
$$
xz^{n}\xi =\varphi(z^{n}) x\xi
$$
from which our claim follows.
\indent{\par}
{\it Case 2.}
Suppose now that $z^{n}$ is in the `middle' of the moment. We need to show that
\begin{equation}
\varphi (w_{1}xz^{n}xw_{2})=\varphi(z^{n})\varphi (w_{1}x^2w_{2})
\end{equation}
for any $w_{1},w_{2}\in {\rm alg}(x,z)$.
Let $h\in {\cal H}_{1}^{0}$ be a vector of norm $\parallel h \parallel =1$.
Using decomposition (8.1) for $b_{1}$ and $b_{2}$, we get
$$
z^{n}h=
\sum_{((i_{1},j_{1}), \ldots , (i_{n},j_{n}))\in I_{2}}
X_{i_{1}}(j_{1})X_{i_{2}}(j_{2})\ldots X_{i_{n}}(j_{n})\,h = b_{2}^{n}h
$$
where $I_{2}$ consists of sequences $((i_{1},j_{1}), \ldots , (i_{n},j_{n}))$
ending with $(i_{n},j_{n})=(2,2)$ and such that for any $2\leq r \leq n$,
either $(i_{r-1},j_{r-1})=(i_{r},j_{r})$, or $|i_{r-1}-i_{r}|=1$ and $j_{r-1}=j_{r}+1$
(this is because $X_{i}(j)X_{i'}(j')=0$ unless $(i,j)=(i',j')$ or $|i-i'|=1$ and $j=j'+1$).
Therefore,
\begin{eqnarray*}
z^{n}h
&=& \langle b_{2}^{n}h, h \rangle \; h \;\;\; {\rm mod}\;({\cal K}^{(2)}(2) \oplus \ldots \oplus {\cal K}^{(n+1)}(2))\\
&=& \varphi (b_{1}^{n})h  \;\;\;\;\;\;\;\,{\rm mod}\;({\cal K}^{(2)}(2) \oplus \ldots \oplus {\cal K}^{(n+1)}(2))
\end{eqnarray*}
where we used Proposition 8.1. This, together with (9.3) and the fact that $z$ and $b_{1}$ are
$\varphi$-identically distributed, implies that
$$
z^{n}(h+\alpha \xi)-\varphi(z^{n})(h+\alpha \xi) \; \perp \; {\mathbb C}\xi \oplus {\cal H}_{1}^{0}
$$
for any $\alpha \in {\mathbb C}$. Since we have $xz^{n}x=pxpz^{n}pxp$
and $p$ is the canonical projection onto $\,{\mathbb C} \xi \oplus {\cal H}_{1}^{0}$,
we get (9.4). This proves that the pair $(x,z)$ is monotone independent w.r.t. $\varphi$.
Since $z$ and $b_{1}$ are $\varphi$-identically distributed,
we get the desired decompositions of the additive free convolution.\hfill $\blacksquare$\\
\indent{\par}
It remains to connect the `orthogonal' decomposition of $\mu \boxright \nu$ given by
Theorem 8.4 with the `boolean' and `monotone' decompositions of $\mu \boxplus \nu$
given by Lemmas 9.1-9.2. We formulate this result, using the associated transforms:
the $K$-transform in the `symmetric' boolean case and the reciprocal Cauchy transform in the
`non-symmetric' monotone case. \\
\indent{\par}
{\sc Theorem 9.3.}
{\it If $\mu$ and $\nu$ are the distributions of $X_{1}$ and $X_{2}$, then}
\begin{eqnarray}
F_{\mu\,\boxplus \,\nu}(z)
&=&
\lim_{m\rightarrow \infty}F_{\mu}(F_{\nu \,\vdash_{m}\mu}(z))\\
K_{\mu \,\boxplus \,\nu}(z)&=&
\lim_{m\rightarrow \infty}(K_{\mu\,\vdash_{m} \nu}(z)+K_{\nu\,\vdash_{m} \mu}(z))
\end{eqnarray}
{\it where the convergence is uniform on compact subsets of ${\mathbb C}^{+}$.}\\
\indent{\par}
{\it Proof.}
For any $m\in {\mathbb N}$, the moments $(\mu_{1}\vartriangleright \mu_{2})(m)$ and
$(\mu_{1}\uplus \mu_{2})(m)$ of compactly supported measures depend only on the moments
of $\mu_{1}$ and $\mu_{2}$ of orders $\leq m$. Therefore, in view of (2.3), Lemma 8.3 and Lemma 9.2,
the moments of $X_1+X_2$ of orders $\leq m$ agree with the corresponding moments of
$\mu\vartriangleright (\nu \vdash_{m-1} \mu)$
(the same holds for the moments of $\nu\vartriangleright (\mu \vdash_{m-1} \nu)$).
Similarly, in view of (2.5), Lemma 8.3 and Lemma 9.1, they agree
with the corresponding moments of $(\mu \vdash_{m-1} \nu) \uplus (\nu \vdash_{m-1} \mu)$.
Therefore,
\begin{eqnarray*}
\mu \boxplus \nu &=& w-\lim_{n\rightarrow \infty} \mu \vartriangleright (\nu \vdash_{n} \mu)\\
\mu \boxplus \nu &=& w-\lim_{n\rightarrow \infty} ((\mu \vdash_{n} \nu)\uplus (\nu \vdash_{n} \mu))
\end{eqnarray*}
which, by Theorem 7.4, gives the assertion.\hfill $\blacksquare$\\
\indent{\par}
{\sc Remark 9.1.}
One can also write (9.5)-(9.6) in the `continued composition' form, similar to (8.4).
In particular, (9.5) gives
\begin{equation}
F_{\mu\,\boxplus \,\nu}(z)
=
F_{\mu}(z-K_{\nu}(z-K_{\mu}(z-K_{\nu}(z-\ldots))))
\end{equation}
which is particularly useful when $\mu$ and $\nu$ have simple $K$-transforms and
computations on continued fractions can be carried out (roughly speaking, the difficulty
in computing the free additive convolution is related to the difficulties arising
in the addition of continued fractions).   \\
\indent{\par}
{\sc Example 9.1.}
We have $\mu \,\boxplus\, \delta_{a}=\mu \vartriangleright (\delta_{a}\,\boxright\, \mu) =
\mu \vartriangleright \delta_{a}\sim \mu_{a}$, where by $\mu_{a}$ we denote the measure
associated with the Jacobi sequences $((\alpha_{0}+a, \alpha_{1}+a, \ldots ), \omega))$.
On the other hand, $\delta_{a} \, \boxplus\, \mu = \delta_{a}\vartriangleright (\mu \,\boxright \,\delta_{a})=\mu_{a}$
since $J(\mu\,\boxright\,\delta_{a})= ((\alpha_{0}, \alpha_{1}+a, \alpha_{2}+a, \ldots ), \omega)$
and $F_{\delta_{a}\vartriangleright \nu}(z)=F_{\nu}(z)-a$ for any $\nu$.\\
\indent{\par}
{\sc Example 9.2.}
Let $F_{\mu}(z)=z-\alpha-\omega/z$ and $F_{\nu}(z)=z-\beta-\gamma/z$. Using Example 8.2 and Theorem 9.3, we obtain
\begin{eqnarray*}
F_{\mu\,\boxplus \,\nu}(z)&=&F_{\mu}(F_{\nu\,\boxright \,\mu}(z))\\
&=&F_{\nu\,\boxright \,\mu}(z)-\alpha-\frac{\omega}{F_{\nu\,\boxright \,\mu}(z)}.
\end{eqnarray*}
This is the algebraic relation which can be used to compute the explicit form of the (reciprocal) Cauchy
transform of $\mu\,\boxplus\,\nu$. Since the solution has a complicated from and it is
known (see, for instance [12]), we shall not give it here.\\
\indent{\par}
{\sc Example 9.3.}
Let $\mu$ be the Wigner measure with mean $\alpha$ and variance $\omega$.
Then it holds that $J(\mu \,\boxplus \,\mu)= ((2\alpha,2\alpha, \ldots ),(2\omega,2\omega, \ldots))$,
which can be shown by using Theorem 9.3 and the result of Example 8.3.
Namely,
\begin{eqnarray*}
[F_{\mu \,\boxplus \,\mu}(z)]_{m}&=&[F_{\mu}(F_{\mu\,\boxright \,\mu}(z))]_{m}\\
&=&
z-J_{m}(z)-\alpha-
\cfrac{\omega}
{z-J_{m-1}(z)-\alpha-
\cfrac{\omega}
{z-J_{m-2}(z)-\alpha-
\cfrac{\omega}{\ldots}}}\\
&=&
z-2J_{m}(z)
\end{eqnarray*}
which proves our assertion.\\
\indent{\par}
{\sc Remark 9.2.}
The decompositions of the free additive convolution of Theorem 9.3 correspond to 
natural decompositions of the free product of graphs. If we consider two rooted
graphs as in Examples 4.1 and 7.1, Theorem 9.3 gives decompositions
\begin{eqnarray*}
{\cal G}_{1}*{\cal G}_{2}&=&{\cal G}_{1}\vartriangleright ({\cal G}_{2}\vdash({\cal G}_{1}\vdash {\cal G}_{2}
\vdash (\ldots )))\\
{\cal G}_{1}*{\cal G}_{2}&=
&\left({\cal G}_{1}\vdash ({\cal G}_{2}\vdash({\cal G}_{1}\vdash (\ldots )))\right)
\star 
\left({\cal G}_{2}\vdash ({\cal G}_{1}\vdash({\cal G}_{2}\vdash (\ldots )))\right)
\end{eqnarray*}
where $\vartriangleright$ and $\star$ denote the so-called comb- and star products of rooted 
graphs, respectively. One can show that both decompositions are related to natural inductive 
definitions of the free product of rooted graphs and to the so-called $m$-free [14]
and $m$-monotone hierarchies [17] of product states.
Recall that the {\it comb product} $({\cal G}_{1}\vartriangleright {\cal G}_{2},e)$ is obtained by attaching 
a replica of ${\cal G}_{2}$ by its root to every vertex of $V_{1}$, whereas the {\it star product}
$({\cal G}_{1}\star {\cal G}_{2},e)$ is obtained by glueing ${\cal G}_{1}$ and ${\cal G}_{2}$
together at their roots. The associated convolutions are the monotone and boolean convolutions,
respectively. Details will be discussed in [1].\\[10pt]

\begin{center}
{\sc Bibliography}
\end{center}
[1] L. Accardi, R. Lenczewski, R. Sa{\l}apata,
Decompositions of the free product of graphs, preprint in preparation, 2006.\\[3pt]
[2] N.I. Akhiezer, {\it The classical moment problem},
Oliver and Boyd, Edinburgh and London, 1965.\\[3pt]
[3] D.~Avitzour, Free products of $C^{*}$- algebras,
{\it Trans.~Amer.~Math.~Soc.} {\bf 271} (1982), 423-465.\\[3pt]
[4] S. Belinschi, Complex Analysis Methods in Noncommutative Probability, Ph.D. Thesis,
Indiana University, 2005.\\[3pt]
[5] H. Bercovici, D. Voiculescu, Free convolution of measuers with unbounded support,
{\it Indiana Univ. Math. J.} {\bf 42} (1993), 733-773.\\[3pt]
[6] Ph. Biane, Processes with free increments, {\it Math. Z.} {\bf 227} (1998), 143-174.\\[3pt]
[7] M. Bo\.{z}ejko, M. Leinert, R. Speicher Convolution and limit theorems
for conditionally free random variables,
{\it Pacific J.Math.} {\bf 175} (1996), 357-388. \\[3pt]
[8] M. Bo\.{z}ejko, A.D. Krystek, {\L}.J. Wojakowski,
Remarks on the {\it r} and $\Delta$ convolution, {\it Math. Z.} {\bf 253} (2006), 177-196\\[3pt]
[9] M. Bo\.{z}ejko, J. Wysocza\'{n}ski, Remarks on $t$-transformations of measures andc convolutions,
{\it Ann.I.H.Poincar\'{e}- PR} {\bf 37}, 6 (2001), 737-761.\\[3pt]
[10] T. Cabanal-Duvillard, V. Ionescu,
Un th\'eor\`eme central limite pour des variables al\'eatoires
non-commutatives, {\it C.\ R.\ Acad.\ Sci.\ Paris},
{\bf 325} (1997), S\'erie I, 1117-1120.\\[3pt]
[11] G.P. Chistyakov, F. Goetze, The arithmetic of distributions in free probability
theory, preprint, arXiv:math.QA/0508245, 2005.\\[5pt]
[12] E. Gutkin, Green's functions of free products of operators with applications to graph
spectra and to random walks, {\it Nagoya Math. J.} {\bf 149} (1998), 93-116.\\[3pt]
[13] Y. Kato, Mixed periodic Jacobi continued fractions, {\it Nagoya Math. J.} {\bf 104} (1986), 129-148.\\[3pt]
[14] R. Lenczewski, Unification of independence in quantum probability,
{\it Infin. Dimens. Anal. Quantum Probab. Relat. Top.}, {\bf 1} (1998), 383-405.\\[3pt]
[15] R. Lenczewski, Reduction of free independence to tensor independence,
{\it Infin. Dimens. Anal. Quantum Probab. Relat. Top.} {\bf 7} (2004), 337-360.\\[3pt]
[16] R. Lenczewski, Noncommutative extension of the Fourier transform and its logarithm,
{\it Studia Math.} {\bf 152} (2002), 69-101. \\[3pt]
[17] R. Lenczewski, R. Sa{\l}apata, Discrete interpolation between monotone probability and free
probability, {\it Infin. Dimens. Anal. Quantum Probab. Relat. Top.} {\bf 9} (2006), 77-106.\\[3pt]
[18] H. Maassen, Addition of freely independent random variables, {\it J. Funct. Anal.} {\bf 106} (1992), 409-438.\\[3pt]
[19] N. Muraki, Monotonic independence, monotonic central limit theorem and monotonic law of small numbers,
{\it Infin. Dimens. Anal. Quantum Probab. Relat. Top.} {\bf 4} (2001), 39-58.\\[3pt]
[20] N. Muraki, Monotonic convolution and monotonic Levy-Hincin formula, preprint,
2001.\\[3pt]
[21] G. Quenell, Combinatorics of free product graphs, {\it Contemp. Math.} {\bf 206}
(1994), 257-281.\\[3pt]
[22] R. Speicher, R. Woroudi, Boolean convolution, in {\it Free Probability
Theory}, Ed. D. Voiculescu, 267-279, {\it Fields Inst. Commun.} Vol.12, AMS, 1997.\\[3pt]
[23] D. Voiculescu, Symmetries of some reduced free product $\mathcal{C}^*$-algebras,
Operator Algebras and their Connections with Topology and Ergodic Theory,
{\it Lecture Notes in Math.} 1132, Springer, Berlin, 1985, 556-588.\\[3pt]
[24] D. Voiculescu, Addition of certain non-commuting random variables,
{\it J. Funct. Anal.} {\bf 66} (1986), 323-246.\\[3pt]
[25] D. Voiculescu, The analogues of entropy and of Fisher's information measure in free probability theory, I,
{\it Commun. Math. Phys.} {\bf 155} (1993), 71-92.\\[3pt]
[26] D. Voiculescu , K. Dykema, A. Nica, {\it Free random variables}, CRM Monograph
Series, No.1, A.M.S., Providence, 1992.
\end{document}